\documentclass[preprint,12pt]{article}
\usepackage[cp1251]{inputenc}
\usepackage[english]{babel}
\usepackage{pgf,tikz,circuitikz}
\usetikzlibrary{arrows}
\usepackage{graphicx}
\usepackage{amssymb}
\usepackage{amsthm}
\usepackage{amsmath}
\usepackage{cmap}
\usepackage{ dsfont }
\usepackage{enumerate}

\usepackage[breaklinks,pdfstartview=FitH,colorlinks,linktocpage,bookmarks=true]{hyperref}
\usepackage{geometry}
\usepackage{pgf,tikz}\usetikzlibrary{arrows}
\newtheorem{theorem}{Theorem}[section]
\newtheorem{lemma}[theorem]{Lemma}

\newtheorem{corollary}[theorem]{Corollary}

\newtheorem{predl}[theorem]{Proposition}

\theoremstyle{definition}

\newtheorem{example}[theorem]{Example}

\theoremstyle{remark}
\newtheorem{remark}[theorem]{Remark}

\renewcommand{\ge}{\geqslant}
\renewcommand{\le}{\leqslant}
\begin{document}
	\title{Transformations of 2-port networks and tiling by rectangles}
	\author{Svetlana Shirokovskikh}
	\date{}
	\maketitle	

\begin{abstract}
	In this paper, we study 2-port networks and introduce new concepts of voltage drop and $\Pi$-equivalence. The main result is that each planar network is $\Pi$-equivalent to a network with no more than 5 edges. This implies that if an octagon in the shape of the letter $\Pi$ can be tiled by squares then it can be tiled by no more than 5 rectangles with rational aspect ratios. Kenyon's theorem from 1998 proves this only for 6 rectangles.
\end{abstract}

\section{Introduction}
\begin{figure}[h]
	\center
	\includegraphics[width=8cm]{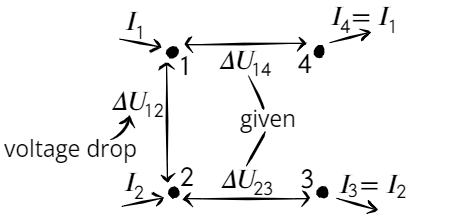}
	\caption{A two-port network (see Section~\ref{prelim2}).}
	\label{fig1}
\end{figure}

We study 2-port networks, that is, weighted graphs with four marked vertices decomposed into two pairs called \emph{ports}. A visual example of such a network is a computer with two USB-ports or a flashlight with two batteries. In two-port circuits, instead of the voltages of marked vertices, we specify the differences of voltages for each port. At the same time, we require that if some current flows into one vertex of a port, then the same current flows out of the second vertex (Figure~\ref{fig1}; we give precise definitions in Section~\ref{prelim2}).

We introduce a new (albeit natural) concept of \emph{voltage drop} between ports and a new \emph{$\Pi$-equivalence} relation that requires not only the same incoming currents in two networks, but also the same voltage drop between ports.

These concepts are interesting from the engineering point of view. If we perform network transformations and keep the incoming currents only, without preserving the voltage drop, then the voltage difference between the vertices from different ports can become large. Informally, a short circuit may occur, disrupting the normal operation of the network. For example, a computer can burn out if you re-solder it, keeping only voltages and currents through the USB-ports, but not the voltage drop. $\Pi$-equivalent transformations avoid such problems.

The main result of this work is that each planar network is $\Pi$-equivalent to a network with no more than 5 edges (Theorem~\ref{equivalence H or I}), and the number 5 is minimal (Example~\ref{example_min=5}). 
Interestingly, if we require preservation of just incoming currents without preservation of the voltage drop, then we can end up with only 3 edges (Remark~\ref{equivalence 3}).

\begin{figure}[h!]
	\includegraphics[width=15cm]{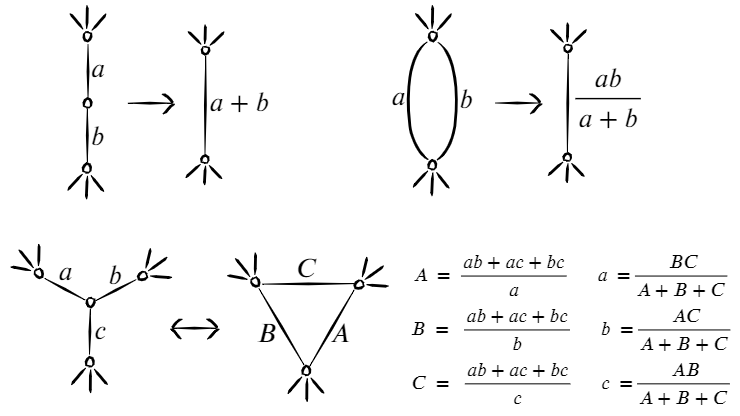}	
	\caption{Elementary transformations (see Section~\ref{prelim1}).}
	\label{elementary}
\end{figure}

\begin{figure}[h]
	\center
	\includegraphics[width=11cm]{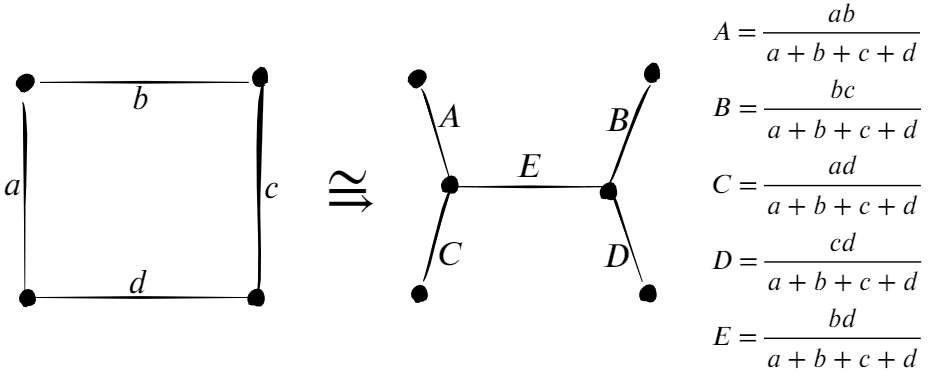}
	\caption{Box-H transformation (see Corollary~\ref{Box-H transformation}).}
	\label{Box-H_transformation}
\end{figure}

The proof uses $\Pi$-equivalent transformations to reduce the number of edges in the network. The simplest of these transformations are well known; they are shown in Figure~\ref{elementary} (see~\cite{CM,K,ZSU}). We introduce a new \textit{Box-H transformation} (Figure~\ref{Box-H_transformation}).

New results on voltage drop and $\Pi$-equivalence make it possible to solve new tiling problems (thanks to the well-known physical interpretation of tilings by R.~L.~Brooks, K.~A.~B.~Smith, A.~G.~Stone and W.~T.~Tutte~\cite{BSST}, described in detail in~\cite{PS}). They imply that if an octagon in the shape of the letter $\Pi$ can be tiled by squares then it can be tiled by no more than 5 rectangles with rational aspect ratios (Theorem~\ref{final_theorem}). The known Kenyon theorem (Theorem~\ref{Kenyon}) proves this only for $6$ rectangles.

\subsection{Transformations of two-port networks}
Transformations of two-port networks that keep incoming currents have been studied in the literature \cite{BIK,B,Bess}. We will give some examples.

\begin{figure}[h]
	\center
	\includegraphics[width=6cm]{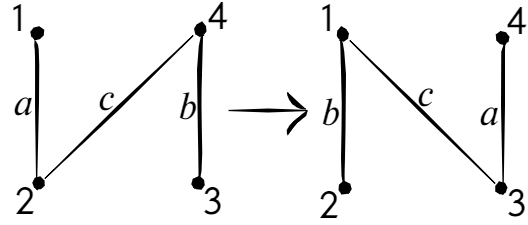}
	\caption{Two-port networks with the same incoming currents.}
	\label{N-N}
\end{figure}

For instance, it is easy to see that for the same differences of incoming voltages the two-port networks in Figure~\ref{N-N} have the same incoming currents. Indeed, if in the network to the right we swap the vertices in each port ($\{1, 4\}$ and $\{2, 3\}$), then we get the network to the left, but differences of incoming voltages change their sign. Therefore, the current flowing out of vertex 1 of the left network is equal to the current flowing into vertex 4 of the right network. In two-port networks, we require that if some current flows into one vertex of a port, then the same flows out of the second vertex. Therefore, in the left network, the current flowing out vertex 1 is equal to the current flowing into vertex 4. By transitivity, the current flowing into vertex 4 is the same in the left and right networks. For the other vertices, the proof is similar.

We say that a two-port network with ports  $\{1, 4\}$ and $\{2, 3\}$ is \emph{symmetric}, if there exists an automorphism of a weighted graph which swap vertices in each pair $\{1, 2\}$ and $\{3, 4\}$.

\begin{figure}[h]
	\center
	\includegraphics[width=11cm]{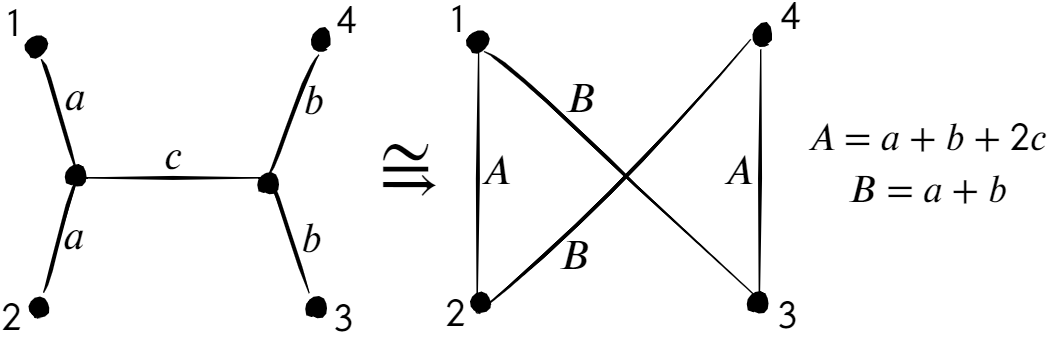}
	\caption{Two-port networks with the same incoming currents.}
	\label{wiki}
\end{figure}

In 1927 Barlett~\cite{B} proved that for any symmetric two-port network there exists a two-port network with four edges ${12}, {34}, {13},$ and ${24}$ with the same incoming currents (Figure~\ref{wiki}).

\subsection{Tilings}

In this section, we give a short overview of known results on tilings by squares.

A natural question is: which polygons can be tiled by squares? The answer for rectangles was obtained by Max Dehn in 1903.

\begin{theorem} [Dehn, 1903\label{dehn}] \textup{\cite{D}} 
	A rectangle can be tiled by squares \textup{(}not necessarily equal\textup{)} if and only
	if the ratio of two orthogonal sides of the rectangle is rational.
\end{theorem}

In this paper, we consider only tilings by finitely many arbitrary squares.

In the case of hexagons the answer was obtained by
R. Kenyon almost a hundred years later.

\begin{theorem} [Kenyon, 1998] \label{th-kenyon}
	\textup {(Cf.~\cite[Theorem~9]{K})} Let $A_1A_2A_3A_4A_5A_6$ be a hexagon with right angles whose vertices are enumerated counterclockwise starting from the vertex of the nonconvex angle. The hexagon $A_1A_2A_3A_4A_5A_6$ can be tiled by squares if and only if the system
	\begin{equation}
		\begin{cases}
			\label{syst}
			A_3A_4\cdot x+ A_1A_2\cdot y=A_2A_3,\\
			A_5A_6\cdot z-A_1A_2\cdot y=A_6A_1;
		\end{cases}
	\end{equation}
	has a nonnegative rational solution $x,y,z$.
\end{theorem}

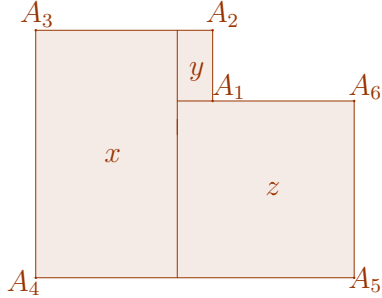
\begin{figure}[htbp]
	\center
	\definecolor{zzttqq}{rgb}{0.6,0.2,0}
\begin{tikzpicture}[line cap=round,line join=round,>=triangle 45,x=0.5cm,y=0.5cm]
\clip(-4.88,-1.69) rectangle (5.72,6.85);
\fill[color=zzttqq,fill=zzttqq,fill opacity=0.1] (1,6) -- (-4,6) -- (-4,-1) -- (5,-1) -- (5,4) -- (1,4) -- cycle;
\draw [color=zzttqq] (1,6)-- (-4,6);
\draw [color=zzttqq] (-4,6)-- (-4,-1);
\draw [color=zzttqq] (-4,-1)-- (5,-1);
\draw [color=zzttqq] (5,-1)-- (5,4);
\draw [color=zzttqq] (5,4)-- (1,4);
\draw [color=zzttqq] (1,4)-- (1,6);
\draw [color=zzttqq] (0,6)-- (0,3.04);
\draw [color=zzttqq] (0,-1)-- (0,3.48);
\draw [color=zzttqq] (0,4)-- (1,4);
\draw[color=zzttqq] (-2.35,2.94) node[anchor=north west] {$x$};
\draw[color=zzttqq] (0.05,5.4) node[anchor=north west] {$y$};
\draw[color=zzttqq] (2.2,2) node[anchor=north west] {$z$};
\fill [color=zzttqq] (1,6) circle (0.5pt);
\draw[color=zzttqq] (1.28,6.45) node {$A_2$};
\fill [color=zzttqq] (-4,6) circle (0.5pt);
\draw[color=zzttqq] (-3.96,6.45) node {$A_3$};
\fill [color=zzttqq] (-4,-1) circle (0.5pt);
\draw[color=zzttqq] (-4.34,-1.14) node {$A_4$};
\fill [color=zzttqq] (5,-1) circle (0.5pt);
\draw[color=zzttqq] (5.32,-1.06) node {$A_5$};
\fill [color=zzttqq] (5,4) circle (0.5pt);
\draw[color=zzttqq] (5.32,4.4) node {$A_6$};
\fill [color=zzttqq] (1,4) circle (0.5pt);
\draw[color=zzttqq] (1.46,4.4) node {$A_1$};
\end{tikzpicture} 
	\caption{\textup{\cite{PS}} The hexagon $A_1A_2A_3A_4A_5A_6$.}
	\label{ell}
\end{figure}
Figure~\ref{ell} explains why if system~(\ref{syst}) has a rational solution $x,y,z\ge 0$ then one can dissect the hexagon into $3$ rectangles of aspect ratios $x,y,z$ and then into squares.

A polygon is called \emph{orthogonal} if its sides are parallel to the coordinate axes.

We say that a polygon is \emph{generic}, if the $x$-coordinates of the sides parallel to the $y$-axis are pairwise distinct.

\begin{theorem} \label{corollary_about_6}
	(Corollary from Theorem~\ref{Kenyon} and Lemma~\ref{correspondence-general} below; Kenyon, 1998)
	Let P be a generic orthogonal polygon with 2n sides. If it can be tiled by squares then it can be tiled by not more than $n(n-1)/2$ rectangles with rational aspect ratios.
	In particular, a generic orthogonal hexagon and octagon can be tiled by not more than 3 and 6 rectangles with rational side ratios, respectively.
\end{theorem}

In this paper we give a proof of the following stronger assertion on a $\Pi$-shaped octagon.

\begin{figure}[h]
	\center
	\includegraphics[width=10cm]{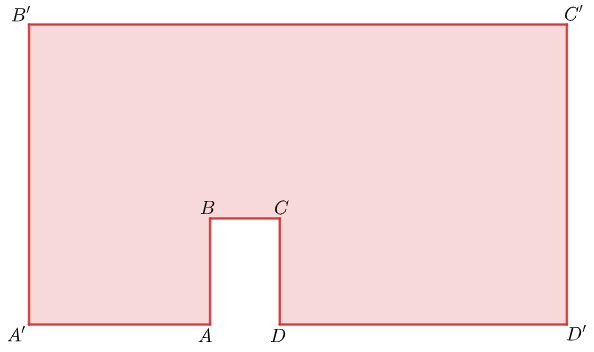}
	\caption{A $\Pi$-shaped octagon.}
	\label{P-8}
\end{figure}

Let $ABCD$ and {$A'B'C'D'$} be rectangles such that $A$ and $D$ lie on {$A'D'$} whereas $B$ and $C$ lie strictly inside {$A'B'C'D'$}. Then the closure of the complement {$A'B'C'D'\setminus ABCD$} is called a \textit{$\Pi$-shaped octagon} (Figure~\ref{P-8}).

\begin{theorem}\label{final_theorem}
	If a $\Pi$-shaped octagon can be tiled by squares, then it can be tiled by no more than 5 rectangles with rational aspect ratios.
\end{theorem}

\subsection{Organization of the paper}
In Section~\ref{prelim}, we recall known definitions and theorems on electrical networks. In Subsection~\ref{prelim1}, we give rigorous definitions of an electrical network, its response, and equivalence. In addition, we introduce some equivalent transformations. In Subsection~\ref{prelim2}, we give a definition of a two-port network and generalize the concept of response to this set up. In short Subsection~\ref{prelim3}, we give a lemma that connects tilings and electrical networks.

In Section~\ref{new_equivalence}, we introduce the new concepts of voltage drop and $\Pi$-equivalence. Then we prove the $\Pi$-equivalence of a Box-network and an H-network (Figure~\ref{Box-H_transformation}). Then we state and prove the main theorem of this paper (Theorem~\ref{equivalence H or I}). After that, we give an example showing that the assumptions of the theorem are sharp.

In Section~\ref{tiling_octagons}, using Theorem~\ref{equivalence H or I}, we prove Theorem~\ref{final_theorem} on tilings of a $\Pi$-shaped octagon. Finally, we give an example of an octagon which cannot be tiled by less than 5 rectangles with rational aspect ratios.

\section{Preliminaries}\label{prelim}
\subsection{Electrical networks}\label{prelim1}
An \textit{electrical network with ${t}$ terminals} is a connected graph with a nonnegative real number (\textit{conductance}) assigned to each edge, and ${t}$ marked (\textit{boundary}) vertices. For simplicity we assume that the graph has neither multiple edges nor loops. Generalizations for graphs with multiple edges are left to the reader.

We say that an electrical network is \textit{planar}, if the graph is embedded in the unit disc in such a way that the boundary vertices belong to the boundary of the disc and the edges have no common points except their common vertices.

Fix an enumeration of the vertices ${1}$, ${2}$, ${\dots}$, ${n}$ such that ${1, \dots, t}$ are the boundary ones. If an electrical network is planar then we assume that the boundary vertices are enumerated counterclockwise along the boundary of the unit disc. Denote by ${kl}$ the edge between the vertices ${k}$ and ${l}$. Denote by $c_{kl}$ the conductance of the edge ${kl}$. Set $c_{kl} = 0$, if there is no edge between ${k}$ and ${l}$.

An \textit{electrical circuit} is an electrical network along with ${t}$ real numbers ${U_1, \dots, U_t}$ (\textit{incoming voltages}) assigned to the boundary vertices.

Each electrical circuit gives rise to certain numbers $U_k$, where $1\le k\le n$ (\textit{voltages} at the vertices), and $I_{kl}$, where $1\le k,l\le n$ (\textit{currents} through the edges, where ${I_{kl}} = 0$ if there is no edge $kl$ according to axiom~\textup{(C)} below).
These numbers are defined by the following $2$ axioms:

\begin{enumerate}[(1)]
	\item[(C)]\textit{The Ohm law.} For each pair of vertices $k,l$ we have $I_{kl}=c_{kl}(U_k-U_l)$.
	\item[(I)]\textit{The Kirchhoff current law.} For each vertex $k>t$ we have $\sum_{l=1}^{n} I_{kl}=0$.
\end{enumerate}

The numbers $U_k$ and $I_{kl}$ are well-defined by these axioms by the following classical result.

\begin{theorem} \textup{~\cite[Theorem~2.1.]{PS}} For any electrical circuit
	the system of linear equations \textup{(C), (I)} in variables $U_k$, where $t < k\le n$, and $I_{kl}$, where $1\le k,l\le n$, has a unique solution.
\end{theorem}

The reciprocal of the conductance is called \emph{the resistance}. Denote by $R_{kl}$ the resistance of the edge ${kl}$.
Numbers $(I_{1}, \dots, I_{t}) := ({{\sum_{k=1}^{n} I_{1k}}}, \dots, {{\sum_{k=1}^{n} I_{tk}}})$ are called \textit{incoming currents}.
The linear map $\mathds{R}^t \to \mathds{R}^t$: $({U_{1}}, \dots, {U_{t}}) \mapsto (I_{1}, \dots, I_{t})$ is called \emph{the response} of the network. The description of the responses of planar networks is very interesting (see~\cite{CM},~\cite{K}).
Two networks are called \textit{equivalent} if their responses are equal.

The results of this section are well-known and easily deduced from the definitions (see~\cite[Theorem on electric transformations in \S2.3]{SSU}).

\begin{figure}[h]
	\center
	\includegraphics[width=8cm]{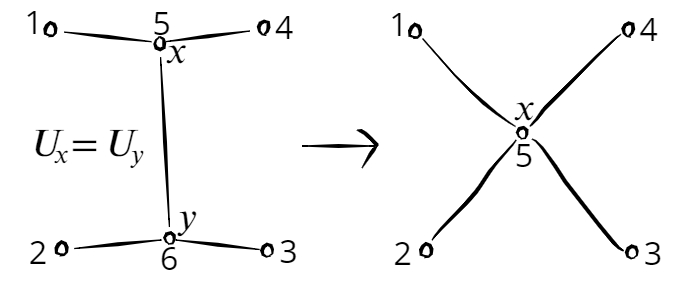}
	\caption{Combining the vertices $x$ and $y$ with the same voltages in an I-network.}
	\label{same_potentials2}
\end{figure}

\begin{predl}\label{same potential0} (Figure~\ref{same_potentials2})
	Let $xy$ be an edge such that $U_x=U_y$ for any incoming voltages. For each $z \neq x$ replace each edge $zy$ in the network by an edge $zx$ of the same resistance and then remove the edge $xy$ and the vertex $y$. Then we get an equivalent network.
\end{predl}

This transformation is called \emph{combining the vertices $x$ and $y$ with the same voltages}.

\begin{remark}\label{remark} (Parallel connection, Figure~\ref{elementary}, top right).
If a network has two multiple edges of resistances $a$ and $b$, then we can replace them by one edge of resistance $ab/(a+b)$ between the same vertices. The resulting network is equivalent to the initial one. Using such parallel connections, we can always transform any network with multiple edges into an equivalent network without multiple edges. The latter network has less edges than the former one. If the resistances of edges in the original network are rational, then the resistances in the new network are also rational. This allows us to consider only networks without multiple edges in what follows.
\end{remark}

\begin{predl}\label{series equivalence}(Figure~\ref{elementary}, top left)
For each $a, b > 0$, the following two networks with two terminals are equivalent:
\begin{enumerate}[(1)]
\item[(1)] The network with three vertices and two edges 13 and 23 of resistances
\begin{align*}
R_{13} = a, R_{23} = b.
\end{align*} 
\item[(2)] The network with two vertices and one edge 12 of resistance
\begin{align*}
R_{12} = a+b.
\end{align*}
\end{enumerate}
\end{predl}

\begin{corollary}\label{series transformation} (Figure~\ref{elementary}, top left)
 Let $xy$ and $yz$ be edges such that $y$ is a non-boundary vertex of degree two. Replace those two edges by an edge $xz$ of resistance $R_{xz}:=R_{xy}+R_{yz}$ and then remove the vertex $y$. Then we get an equivalent network.
\end{corollary}
	
This transformation is called \emph{the series connection} applied to the edges $xy$ and $yz$.

\begin{predl}\label{Delta-Y equivalence} (Figure~\ref{elementary}, bottom)
For each $a, b, c > 0$ , the following two networks with three terminals are equivalent:
\begin{enumerate}[(1)]
\item[($\Delta$)] The network with three vertices and three edges 12, 23, 13 of resistances
\begin{align*}
R_{12} = a, R_{23} = b, R_{13} = c.
\end{align*}
		
\item[($Y$)] The network with four vertices and three edges 14, 24, 34 of resistances
\begin{align*}
R_{14} = ac/(a+b+c), R_{24} = ab/(a+b+c), R_{34} = bc/(a+b+c).
\end{align*}
\end{enumerate}
\end{predl}

\begin{corollary}\label{Delta-Y transformation} (Figure~\ref{elementary}, bottom)
Let x, y, and z be vertices such that each pair of these vertices is joined by an edge. Add a new non-boundary vertex $w$ and replace edges $xy$, $xz$, $yz$ by edges $xw, yw, zw$ of resistances
\begin{align*} &R_{xw}:=R_{xy}R_{xz}/(R_{xy}+R_{xz}+R_{yz}),\\
&R_{yw}:=R_{xy}R_{yz}/(R_{xy}+R_{xz}+R_{yz}),\\
&R_{zw}:=R_{xz}R_{yz}/(R_{xy}+R_{xz}+R_{yz}).
\end{align*}
Then we get an equivalent network.
\end{corollary}

This transformation is called \emph{the $\Delta-Y$ transformation} applied to the edges $xy, yz, xz$. The inverse transformation is called \emph{the $Y-\Delta$ transformation}.

 Parallel connection, series connection, $\Delta-Y$ transformation, and $Y-\Delta$ transformation are called \emph{elementary transformations} (Figure.~\ref{elementary}).

\begin{theorem}(Kenyon, 1998) \label{Kenyon} 
	\textup {~\cite[Theorem~8]{K}}
	For any planar network $G$ with $t$ terminals, there exists an equivalent planar network which has no more than $t(t-1)/2$ edges and is obtained from the network $G$ by a sequence of elementary transformations.
\end{theorem}

\subsection{Two-port networks}\label{prelim2}

A \textit{two-port network} is an electrical network with 4 terminals decomposed into 2 pairs called \emph{ports}.

Fix an enumeration of the vertices ${1}$, ${2}$, ${\dots}$, ${n}$ such that ${1}$, ${2}$, ${3}$, ${4}$ are the boundary ones decomposed into the pairs $\{1,4\}$ and $\{2,3\}$.
As before, denote by $c_{kl}$ the conductance of the edge ${kl}$. Set $c_{kl} = 0$, if there is no edge between ${k}$ and ${l}$. A planar two-port network is defined similarly to a planar electrical network. If a two-port network is planar then we assume that the boundary vertices are enumerated counterclockwise along the boundary of the unit disc.

A \textit{two-port circuit} is a two-port network along with ${2}$ real numbers $\Delta{U_{14}}$ and $\Delta{U_{23}}$ (\textit{differences of incoming voltages}), assigned to the pairs of boundary vertices.

Each two-port circuit gives rise to certain numbers $U_k$, where $1\le k\le n$ (\textit{voltages} at the vertices), and $I_{kl}$, where $1\le k,l\le n$ (\textit{currents} through the edges, where ${I_{kl}} = 0$ if there is no edge $kl$).
These numbers are defined by the following $4$ axioms:

\begin{enumerate}[(1)]
\item[(C)]\textit{The Ohm law.} For each pair of vertices $k,l$ we have $I_{kl}=c_{kl}(U_k-U_l)$.
\item[(I)]\textit{The Kirchhoff current law.} For each vertex $k>4$ we have $\sum_{l=1}^{n} I_{kl}=0$.
\item[(P)] \textit{Port-isolation condition.} ${{\sum_{l=1}^{n} I_{1l}+\sum_{l=1}^{n} I_{4l}=0}}$.
\item[(B)]\textit{Boundary conditions.} $ U_1 - U_4 = \Delta{U_{14}}$ and $U_2 - U_3 = \Delta{U_{23}}$.
\end{enumerate}

\begin{theorem}\label{exist_uniq}
For any two-port electrical circuit there exist a unique collection ${I_{kl}}$, ${1\le k,l\le n}$, and a unique up to adding a constant collection ${U_{k}}$, ${1\le k\le n}$, that satisfy conditions (C), (I), (P), (B).
\end{theorem}

The proof of Theorem~\ref{exist_uniq} is well-known, but we present it for completeness.
It is obtained similarly to the proof of Theorem~2.1 from \textup{~\cite{PS}} and uses the following lemma.

\begin{lemma} \textup{~\cite[Lemma~5.1]{PS}} \label{energy} \label{fullpower} 
	Consider an electrical network with the vertices  ${1, \dots, n}$ such that ${1, \dots, t}$ are the boundary ones.
	Suppose that the numbers ${U_k}$, where ${1\le k\le n}$, and ${I_{kl}}$, where ${1\le k,l\le n}$, satisfy the Ohm law~\textup{(C)} and the Kirchhoff current law~\textup{(I)} from \S\ref{prelim}. Set ${I_u=\sum_{k=1}^n I_{uk}}$. Then
	$$\sum\limits_{1\le k<l\le n}(U_k-U_l)I_{kl}
	= \sum\limits_{1\le u\le t}U_u I_u.$$
\end{lemma}

\begin{proof}[Proof of Theorem~\ref{exist_uniq}.] 
	
	Take $C \in \mathds{R}$. It is easy to see that a collection of currents ${I_{kl}}$ and a collection of voltages ${U_{k}}$ satisfies properties (C), (I), (P), (B) if and only if the collection of currents ${I_{kl}}$ and the collection of voltages ${U_{k}} + C$ satisfies the same properties. We prove the existence and uniqueness of solution under the additional condition $U_1=0$. 

\textit{Uniqueness}.
Suppose there are two collections of currents ${I^{I,II}_{kl}}$ and
voltages ${U^{I,II}_k}$ satisfying laws (C), (I), (P), (B) and ${U^{I}_1 = U^{II}_2 = 0}$ for some differences of incoming voltages $\Delta{U_{14}}^{I}$ = $\Delta{U_{14}}^{II}$ and $\Delta{U_{23}}^{I}$ = $\Delta{U_{23}}^{II}$. Then their differences ${I_{kl}=I^{I}_{kl}-I^{II}_{kl}}$ and ${U_{k}=U^{I}_{k}-U^{II}_{k}}$ satisfy laws (C), (I), (P), (B) for zero incoming voltages $\Delta{U_{14}}=0$ and $\Delta{U_{23}}=0$.
 Then $U_1 = U_4$ and $U_2 = U_3$. Using $(I)$ and $(P)$, we obtain $I_1=-I_4$ and $I_2=-I_3$. It follows that $\sum_{u=1}^4 U_u I_u = 0$.

Then by Lemma~\ref{energy} we have $$\sum\limits_{1\le k<l\le n}(U_k-U_l)^2 c_{kl}
=\sum\limits_{1\le k<l\le n}(U_k-U_l)I_{kl}
= \sum\limits_{1\le u\le 4}U_u I_u = 0.$$

For each ${k,l}$ we have either ${c_{kl}>0}$ or ${c_{kl}=0}$. Thus each summand $ c_{kl}(U_k-U_l)^2=0$. Since the network is connected it follows that all the voltages ${U_{k}}$ are equal to each other.
But $U_1 = U_1^{I} - U_1^{II}=0$. Hence ${U_k=0}$, ${I_{kl}=0}$, and thus ${I^{I}_{kl}=I^{II}_{kl}}$, ${U^{I}_{k}=U^{II}_{k}}$ for each ${k,l}$.

\textit{Existence}. The number of equations in the system (C), (I), (P), (B)  equals the number of variables (if we fix $U_1=0$).
We have proved that the system has a unique solution for  $\Delta{U_{14}}=0$, $\Delta{U_{23}}=0$. By the finite-dimensional Fredholm alternative it has a solution for each $\Delta{U_{14}}$, $\Delta{U_{23}}$.
\end{proof}

The numbers $(I_{1}, I_{2}) := ({{\sum_{k=1}^{n} I_{1k}}}, {{\sum_{k=1}^{n} I_{2k}}})$ are called \emph{the incoming currents}.
The linear map $\mathds{R}^2 \to \mathds{R}^2$: $(\Delta{U_{14}}, \Delta{U_{23}}) \mapsto (I_{1}, I_{2})$ is called \textit{the response} of the two-port network.

\begin{example}\label{Box_response}
	The two-port network with four vertices and four edges $12, 23, 34, 14$ of resistances
	$R_{12} = a, R_{14} = b, R_{34} = c, R_{23} = d$ (\textit{Box-network}, see Figure~\ref{Box-H_transformation} to the left) has the following response matrix:
	\begin{center}
		${\left( \begin{matrix}
				(a+b+c)/(ab+ac) & -1/(a+c)\\
				-1/(a+c)  & (a+c+d)/(ad+cd)
			\end{matrix}\right)}$.
	\end{center}
\end{example}

\begin{proof}
	Using  $(C), (I), (P), (B),$ we obtain the system of equations:
	
	\begin{equation}
		\label{syst2}
		\begin{cases}
		U_1-U_4=\Delta U_{14}, \\
		U_2-U_3=\Delta U_{23},\\
		I_{14}b=U_1-U_4,\\
		I_{12}a=U_1-U_2,\\
		I_{43}c=U_4-U_3,\\
		I_{23}d=U_2-U_3,\\ 
		I_{12}+I_{14}+I_{23}+I_{43}=0.
		\end{cases}
	\end{equation}
	From this system we have:
	\begin{align*}
		&I_{14}=\Delta{U_{14}}/b,\\
		&I_{12}=(\Delta{U_{14}}-\Delta{U_{23}})/(a+c),\\
		&I_{43}=(\Delta{U_{23}}-\Delta{U_{14}})/(a+c),\\
		&I_{23}=\Delta{U_{23}}/d,\\
		&I_{1} = I_{14}+I_{12}=((a+b+c)\Delta{U_{14}}-b\Delta{U_{23}})/(b(a+c)),\\
		&I_{2} = I_{21}+I_{23}=((a+c+d)\Delta{U_{23}}-d\Delta{U_{14}})/(d(a+c)).
	\end{align*}
	Hence
	\begin{center}
		${\left( \begin{matrix}
				(a+b+c)/(ab+ac) & -1/(a+c)\\
				-1/(a+c)  & (a+c+d)/(ad+cd)
			\end{matrix}\right)} \cdot
		{\left( \begin{matrix}
				\Delta{U_{14}}\\
				\Delta{U_{23}}
			\end{matrix}\right)} = 
		{\left( \begin{matrix}
				{I_{1}}\\
				{I_{2}}
			\end{matrix}\right)}
		$.
		
	\end{center}
	
\end{proof}

\subsection{Electrical networks and tilings}\label{prelim3}

\begin{figure}[h]
	\center
	\includegraphics[width=12cm]{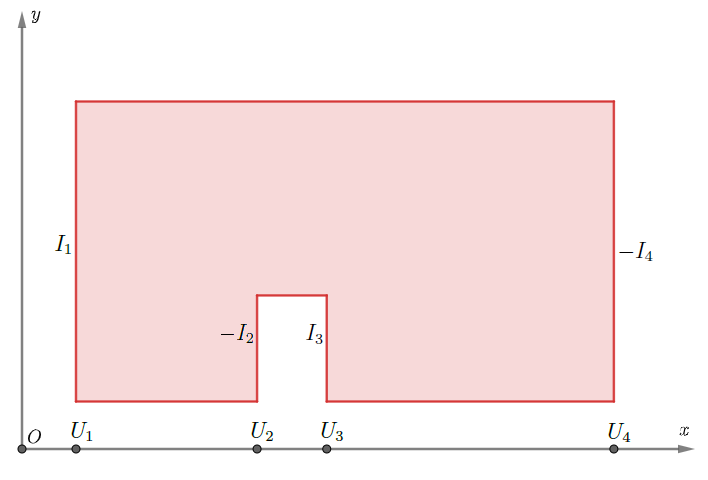}
	\caption{An orthogonal polygon.}
	\label{example_ort}
\end{figure}

Hereafter $P$ is an \emph{orthogonal} polygon (Figure~\ref{example_ort}), i.e., a polygon with the sides parallel to the coordinate axes. Hereafter $P$  is \emph{simple}, i.e., the boundary $\partial P$ is connected. Enumerate the sides of the polygon parallel to the $y$ axis counterclockwise along $\partial P$.
Denote by $t$ the number of these sides.
Let $I_u$ be the signed length of the side $u$, where the sign of $I_u$ is ``$+$'' (respectively, ``$-$'') if the $P$ locally
lies to the right (respectively, left) of the side $u$. Let $U_u$ be the $x$-coordinate of the side $u$. Assume that $P$ is
\emph{generic} in the sense that the numbers ${U_1, \dots, U_t}$ are pairwise distinct.

We say that an edge $kl$ of a circuit is \emph{essential}, if $I_{kl} \neq 0$.

\begin{lemma}\textup{~\cite[Lemma~5.4.]{PS}} \label{correspondence-general}  Let $P$ be a generic simple orthogonal polygon with $t$ vertical sides of
	signed lengths ${I_1, \dots, I_t}$ and ${x\mbox{-co}}$ordinates ${U_1, \dots, U_t}$. Then the following 2 conditions are
	equivalent:
	\begin{enumerate}[(1)]
		\item \label{c-g-1} the polygon $P$ can be tiled by $m$ rectangles with ratios of the horizontal side to the vertical equal to $R_1, \dots, R_m$\textup{;}
		\item \label{c-g-2} there is a planar electrical circuit with $t$ boundary vertices, $m$ essential edges of resistances ${R_1,\dots,R_m}$, incoming voltages ${U_1, \dots, U_t}$, and incoming currents ${I_1, \dots, I_t}$.
	\end{enumerate}
\end{lemma}

\section{$\Pi$-equivalence of two-port networks}\label{new_equivalence}
Now we introduce a new concept, which appears in the paper for the first time.

The linear map $\mathds{R}^2 \to \mathds{R}$: $(\Delta{U_{14}}, \Delta{U_{23}}) \mapsto U_{1}-U_{2}$ is called \textit{the voltage drop} of a two-port network. We say that two-port networks are \textit{$\Pi$-equivalent} 
if their responses and voltage drops are equal. In Lemma~\ref{correspondence-general2} we are going to see how these concepts naturally arise in the study of tilings by rectangles.

\begin{example}\label{Box_drop}
The Box-network (see Example~\ref{Box_response} and Figure~\ref{Box-H_transformation}) has the voltage drop $U_{1}-U_{2} = (\Delta{U_{14}}-\Delta{U_{23}})a/(a+c)$.
\end{example}

\begin{proof}
This follows directly from~(\ref{syst2}).
\end{proof}
We note the following two obvious propositions.

\begin{predl}\label{same potential} (Figure~\ref{same_potentials2}, cf. Proposition~\ref{same potential0})
	Let $xy$ be an edge of a two-port network such that $U_x=U_y$ for any differences of incoming voltages. For each $z \neq x$ replace each edge $zy$ in the network by an edge $zx$ of the same resistance and then remove the edge $xy$ and the vertex $y$. Then we get an $\Pi$-equivalent network.
\end{predl}

This transformation, which generalizes the combining of vertices with the same voltages, is still called \emph{combining the vertices $x$ and $y$ with the same voltages}.

\begin{predl}\label{equel=>}
	If two electrical networks with four terminals are equivalent, then the two-port networks with the same edges and the same resistances are $\Pi$-equivalent.
In particular, elementary transformations preserve the $\Pi$-equivalence classes of two-port networks.
\end{predl}

	Our main tool is a new transformation of two-port networks that preserves both the response and the voltage drop (Figure~\ref{Box-H_transformation}). This transformation in particular shows that the reciprocal of Proposition~\ref{equel=>} is not true (one can check directly by computing the response matrices that it is not an equivalent transformation).
	A network with six vertices and five edges ${15}, {25}, {46}, {36}, {56}$ is called \textit{an H-network}.
	
\begin{theorem}\label{Box-H}(Figure~\ref{Box-H_transformation})
For each $a, b, c, d > 0$ the following two-port networks are equivalent:
\begin{enumerate}[(1)]
\item[(Box)] 
The network with four vertices and four edges $12, 23, 34, 14$ of resistances
\begin{align}
\nonumber
R_{12} = a, R_{14} = b, R_{34} = c, R_{23} = d.	\end{align}
		
\item[(H)] The network with six vertices and five edges $15, 25, 36, 46, 56$ of resistances
\begin{align*}
& R_{15} = ab/(a+b+c+d), R_{25} = ad/(a+b+c+d), R_{46} = bc/(a+b+c+d),\\
& R_{36} = cd/(a+b+c+d), R_{56} = bd/(a+b+c+d). \end{align*}
\end{enumerate}
\end{theorem}

\begin{proof}
From Examples~\ref{Box_response} and \ref{Box_drop} we know the response and the voltage drop for the Box-network. Using  $(C), (I), (P), (B),$ we obtain the system of equations for the H-network:
\begin{equation*}
	\begin{cases}
U_1-U_4=\Delta U_{14},\\
U_2-U_3=\Delta U_{23},\\
I_{15}ab/(a+b+c+d)=U_1-U_5,\\
I_{46}bc/(a+b+c+d)=U_4-U_6,\\
I_{56}bd/(a+b+c+d)=U_5-U_6,\\
I_{25}ad/(a+b+c+d)=U_2-U_5,\\
I_{36}cd/(a+b+c+d)=U_3-U_6,\\
I_{15}+I_{46}=0,\\
I_{15}+I_{25}-I_{56}=0,\\
I_{46}+I_{56}+I_{36}=0.
\end{cases}
\end{equation*}
Using system of equations, we get
\begin{align*}
&I_{15}=I_{64}=((a+b+c)\Delta{U_{14}}-b\Delta{U_{23}})/(b(a+c)),\\
&I_{25}=I_{63}=((a+c+d)\Delta{U_{23}}-d\Delta{U_{14}})/(d(a+c)),\\
&U_{1}-U_{2}=(\Delta{U_{14}}-\Delta{U_{23}})a/(a+c).
\end{align*}

The voltage drops and responses of the networks are the same, thus the networks are $\Pi$-equivalent.
\end{proof}	

The following corollaries are obtained similarly to \cite[Proof of Theorem on electrical transformations, page 35]{SSU}.
\begin{corollary}\label{Box-H transformation1} (Figure~\ref{case2new})
Let $xy, yz, zt, tx$ be edges of a two-port network such that $I_{xy}+I_{tz}=0$ for any differences of incoming voltages.
Add 2 new non-boundary vertices $v$ and $w$ and replace edges $xy, yz, zt, tx$ by edges $xv, yv, tw, zw, vw$ of resistances
	
\begin{align}\label{eq1} 
	& R_{xv} = R_{xy}R_{xt}/(R_{xy}+R_{xt}+R_{tz}+R_{yz}),\\
	\nonumber
	&R_{yv} = R_{xy}R_{yz}/(R_{xy}+R_{xt}+R_{tz}+R_{yz}), \\
	\nonumber
	&R_{tw} = R_{xt}R_{tz}/(R_{xy}+R_{xt}+R_{tz}+R_{yz}),\\
	\nonumber
	& R_{zw} = R_{tz}R_{yz}/(R_{xy}+R_{xt}+R_{tz}+R_{yz}),\\ 
	\nonumber
	&R_{vw} = R_{xt}R_{yz}/(R_{xy}+R_{xt}+R_{tz}+R_{yz}). \end{align}
Then we get a $\Pi$-equivalent network.
\end{corollary}

\begin{corollary}\label{Box-H transformation} (See the first transformation in Figure~\ref{case2new})
Let $xy, yz, zt, tx$ be edges of a two-port network such that 
if we remove $xy, yz, zt, tx$ then any two vertices, where the first one is from the set $\{1, 4, x, t\}$ and the second one from the set $\{2, 3, y, z\}$, are in the different connected components.
Add 2 new non-boundary vertices $v$ and $w$ and replace edges $xy, yz, zt, tx$ by edges $xv, yv, tw, zw, vw$ of resistances given by formula~\eqref{eq1}. Then we get a $\Pi$-equivalent network.
\end{corollary}

This transformation is called \emph{the Box-H transformation} applied to the edges $xy, yz, zt, tx$.

\begin{theorem}\label{equivalence H or I}
	For each planar two-port network $G$, there exists a $\Pi$-equivalent planar two-port network $G'$ satisfying the following two conditions:

\begin{enumerate}[(1)]
	\item[(1)] The network $G'$ has no more than 4 edges or it is the H-network.
	
	\item[(2)] The network $G'$ can be obtained from the network G by elementary transformations (Figure~\ref{elementary}), combining vertices with the same voltages (Figure~\ref{same_potentials2} and Proposition~\ref{same potential}), and Box-H transformations (Figure~\ref{Box-H_transformation}).
\end{enumerate}
\end{theorem}

\begin{remark}\label{rational}
	If all the edges of the network $G$ have rational resistances, then the resistances of the edges of the network $G'$ are also rational.
\end{remark}

If we require preservation of just network response without preservation of the voltage drop, then we can end up with only 3 edges:

\begin{remark}\label{equivalence 3}
	For each planar two-port network, there exists a planar two-port network with the same response which has no more than 3 edges. (One can prove this using Theorem~\ref{equivalence H or I} and equality of the responses of the networks in Figure~\ref{N-N}.)
\end{remark}

\begin{proof}[Proof of Theorem~\ref{equivalence H or I}.]
	Consider all planar two-port networks with the minimal number of edges that can be obtained from the network $G$ by a sequence of transformations from Theorem~\ref{equivalence H or I}. Among such networks, choose the network with the maximal number of vertices. Then the obtained \emph{minimal} network has no cycles of length 3, otherwise the number of vertices can be increased using the $\Delta-Y$ transformation. In addition, there are no non-boundary vertices of degree two, otherwise we can reduce the number of edges by a series connection.
	
	By Theorem~\ref{Kenyon} and Proposition~\ref{equel=>}, this network has no more than 6 edges. If it has no more than 4 edges, then the theorem is proved. Otherwise, this network has 5 or 6 edges. Consider those two cases separately.

\textit{Case 1: the network has 5 edges.} This network is either a tree or it contains cycles of length more than three.

\textit{Subcase 1: the network has no cycles.} Then the network has 6 vertices. The vertices 5 and 6 are non-boundary vertices of degree no less than 3. There are 5 edges in total, hence there is the edge 56. Then this is either an H-network or \emph{I-network} (a network with 6 vertices and 5 edges ${15}, {45}, {26}, {36}, {56},$ see Figure~\ref{same_potentials2}). If this is an I-network, then, by the Port isolation (condition (P) in \S\ref{prelim2}) we have $I_{15}=I_{54}$, and hence, by the Kirchhoff current law (condition (I)) we get $I_{56} = 0$. Hence by the Ohm law (condition (C)) we have $U_5 = U_6$, thus the number of edges can be reduced by combining the vertices 5 and 6 with the same voltages (Proposition~\ref{same potential}).
			
\textit{Subcase 2: the network has a cycle.} There is no cycle of length 5, because otherwise there is a non-boundary vertex of degree two. Then there is a cycle of length 4 and one edge with the endpoint outside the cycle (Figure~\ref{case2new}). The non-boundary vertex has degree 3 (otherwise the number of edges can be reduced). Without loss of generality, the cycle does not include the vertex $4$. Performing a Box-H transformation (Theorem~\ref{Box-H}) and then a series connection (Proposition~\ref{series transformation}) we get an H-network (Figure~\ref{case2new}).

\begin{figure}[h]
	\center
	\includegraphics[width=11cm]{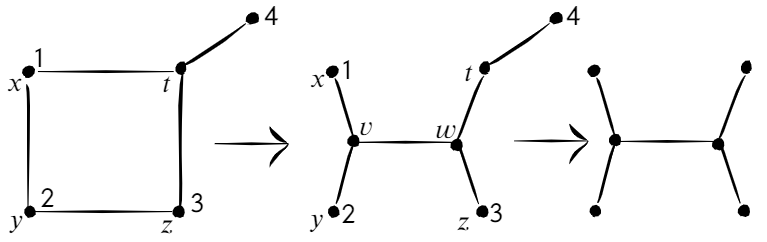}
	\caption{Subcase 2: a network with 5 edges and a cycle of length 4.}
	\label{case2new}
\end{figure}

\textit{Case 2: the network has 6 edges.}
By Lemma~\ref{>=2 edges} (proved right after the proof of the theorem), each boundary vertex has degree at least 2. Each non-boundary vertex has degree at least 3. Hence there are no vertices
of degree 1, thus this network is not a tree. Then there is a cycle of length at least 4. 

There is no cycle of length 6 because non-boundary vertices have degree at least 3.
There is no cycle of length 5 because otherwise the sixth edge either forms a cycle of length 3 or its endpoint has degree 1.

Thus we have a cycle of length 4 with edges, say, $xy$, $yz$, $zt$, $tx$. Then there is also an edge, say, $xw$ with the endpoint outside the cycle. There are no vertices of degree 1 and no cycles of length 3, thus there is an edge $wz$ (Figure~\ref{xyztw}). Then $w$, $y$, and $t$ are boundary vertices because their degree is two. But then the network is not planar. This contradiction proves the theorem (modulo the following lemma).
\end{proof}
\begin{figure}[h]
	\center
	\includegraphics[width=4cm]{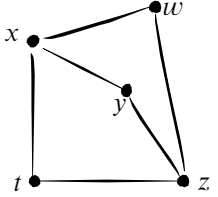}
	\caption{Case 2: a network with 6 edges and a cycle of length 4.}
	\label{xyztw}
\end{figure}

\begin{lemma}\label{>=2 edges}
	Let a two-port network with 6 edges have a boundary vertex of degree 1. Then we can reduce the number of edges by a sequence of transformations from Theorem~\ref{equivalence H or I}.
\end{lemma}

\begin{proof}
	Without loss of generality, let vertex 1 have degree 1. Then the network has either the edge 12, or 13, or 14, or an edge with the endpoint 1 and another non-boundary endpoint (we can denote the latter endpoint by 5).
	
	\textit{Case 1: the edge 12 (Figure~\ref{12}).} 
	Consider the electrical network with 4 terminals with the same edges of the same resistances as the two-port network. Remove the edge 12 and the vertex 1. We get a network with 3 terminals and 5 edges. By Theorem~\ref{Kenyon} it can be transformed into a network with no more than 3 edges by a sequence of elementary transformations. If we return the edge 12 and the vertex 1, then the same transformations remain equivalent ($I_{12}$ depends only on the voltages $U_1$ and $U_2$, and the remaining currents depend on $U_2$, $U_3$, $U_4$ only). Under such transformations, the $\Pi$-equivalence class of the two-port network is preserved by Proposition~\ref{equel=>}, and the number of edges is reduced.
	
	\begin{figure}[h!]
		\center
		\includegraphics[width=8cm]{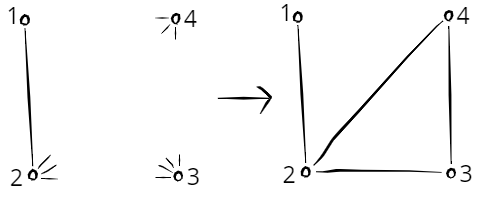}
		\caption{Case 1: the edge 12.}
		\label{12}
	\end{figure}
	
	\textit{Case 2: the edge 14.} This is similar to Case 1.
	
	\textit{Case 3: the edge 13 (Figure~\ref{13}).} 
	The boundary vertices of any planar two-port network belong to the boundary of the unit disc. The edge 13 divides the disk into 2 parts (right and left ones). Consider the electrical network with two boundary vertices 3 and 4, having the same edges and vertices as the edges and vertices to the right of the edge 13 in the two-port network. By Theorem~\ref{Kenyon}, we can transform the resulting electrical network into a network with one edge. Do the same with the left part. We obtain a network with no more than 3 edges.
	
	\begin{figure}[h!]
		\center
		\includegraphics[width=7cm]{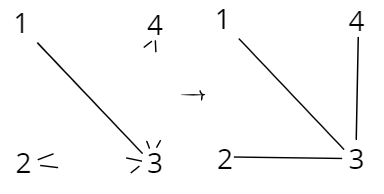}
		\caption{Case 3: edge 13.}
		\label{13}
	\end{figure}
	
	\textit{Case 4: the edge 15 (Figure~\ref{15}).} 
	Remove the edge 15. By Case 1 of Theorem~\ref{equivalence H or I}, which has already been proved, the remaining part of the network can be $\Pi$-equivalently transformed into a network with no more than four edges or into an H-network. In the latter case, we perform a series connection with edge 15 and reduce the number of edges.
\end{proof}
\begin{figure}[h!]
	\center
	\includegraphics[width=6cm]{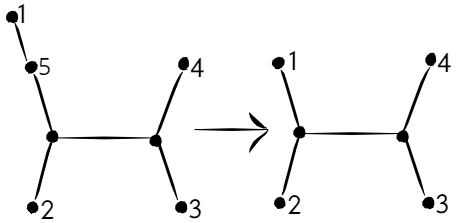}
	\caption{Case 4: edge 15.}
	\label{15}
\end{figure}

Theorem~\ref{equivalence H or I} is proved, but another question arises:
are the conclusions of the theorem sharp? Can  we replace condition (1) by \textit{``the network $G'$ has no more than 4 edges''}? Let us show that the theorem does not remain true in such a stronger form. Moreover, we show that we cannot reduce the number of edges in the two-port network to 4 by any transformations that preserve the response, voltage drop and ``rationality of resistances''.

The linear map $\mathds{R}^2 \to \mathds{R}$: $(\Delta{U_{14}}, \Delta{U_{23}}) \mapsto U_{3}-U_{4}$ is called \textit{the second voltage drop}. For $\Pi$-equivalent networks, the second voltage drops are equal since $U_{3}-U_{4} = \Delta{U_{14}} - \Delta{U_{23}} - (U_{1}-U_{2} ).$ The differences $U_1-U_2$ and $U_3-U_4$ themselves are called \textit{the voltage drop value} and \textit{the second voltage drop value}.

\begin{predl}\label{contrexample} 
	Consider an H-circuit (that is, an H-network with fixed differences of incoming voltages) with edges of rational resistances. Let the ratios of any two numbers among the voltage drop value ($U_1-U_2$), the second voltage drop value ($U_3-U_4$), the incoming currents ($I_1$ and $I_2$) and the differences of incoming voltages ($\Delta{U_{14} }$ and $\Delta{U_{23}}$) are irrational.

Then any planar two-port circuit with edges of rational resistances, with the same differences of incoming voltages, the same voltage drop, and the same response as the H-circuit, has no less than five edges.
\end{predl}

\begin{proof}
	Assume the converse: there exists a two-port circuit $H'$ (and the corresponding network) with less than 5 edges. Consider a two-port network with the minimal number of edges that can be obtained from the two-port network $H'$ by a sequence of transformations from Theorem~\ref{equivalence H or I}. From such networks, choose the network with the maximal number of vertices. Then it has three or four edges of rational resistances. Since for $\Pi$-equivalent networks the voltage drop, the second voltage drop and the response are the same, then in the resulting network all ratios from the assumption of the proposition are also irrational. Let us consider two cases: the resulting network has no cycles or it has a cycle. 

\textit{Case 1: the network has no cycles.} Then the network is a tree with 4 or 5 vertices. Since a non-boundary vertex cannot have degree one, it follows that some boundary vertex has degree one. Let, say, vertex 1 have degree one. There is neither edge 12 nor 14, because the resistances of edges are rational, whereas $(U_1-U_2)/I_1$ and $(U_1-U_4)/I_1$ are irrational. Then there is either the edge 13 or 15.

\textit{Subcase 1.1: the edge 13.} The edge 13 divides the planar network into two parts. Without loss of generality, there are no non-boundary vertices in the part containing the boundary vertex 4. This part of the network can only have edges 14 or 34, because the network is planar. Since the network has the edge 13 and no cycles, vertex 4 has degree 1, and there is either the edge 14 or 34 of rational resistance. But $(U_1-U_4)/I_1$ and $(U_3-U_4)/I_1$ are irrational. This is a contradiction.

\textit{Subcase 1.2: the edge 15.} 
The non-boundary vertex 5 has degree at least 3. Without loss of generality, there are the edges 15, 25, 35. Then the vertex 4 has degree 1, and there is neither edge 14 nor 34, because the resistances of edges are rational, whereas $(U_1-U_4)/I_1$ and $(U_3-U_4)/I_1$ are irrational. Then there is the edge 45 (Figure~\ref{X}). By the rationality of the resistances of the edges 15 and 45 and the port isolation we get a contradiction with the irrationality of the ratio $\Delta{U_{14}}/I_1$.

\begin{figure}[h]
	\center
	\includegraphics[width=4cm]{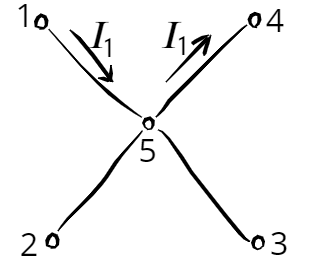}
	\caption{Subcase 1.2.}
	\label{X}
\end{figure}

\textit{Case 2: the network has a cycle.} Since there are no cycles of length 3, this is a Box-network. The rational resistance of the edge 12 is equal to $(U_1-U_2)/(I_1-I_{14})$, the rational resistance of the edge 34 is equal to $(U_3-U_4)/(I_1-I_{14})$. We get a contradiction with the irrationality of the ratio $(U_1-U_2)/(U_3-U_4)$.
\end{proof}

\begin{example}\label{example_min=5}
Consider the H-network with the edge resistances $R_{15}=4, R_{25}=1, R_{46}=4, R_{36}=2, R_{56}=2$. Let differences of incoming voltages be $\Delta{U_{14}} = 10 - 2 \sqrt{2}$ and $\Delta{U_{14}} = 8 - 5 \sqrt{2}$. Then the assumptions of Proposition~\ref{contrexample} are satisfied. Thus each network obtained from this network by the transformations from Theorem~\ref{equivalence H or I} has at least 5 edges.
\end{example}

\section{Tilings of octagons}\label{tiling_octagons}

\begin{figure}[h]
	\center
	\includegraphics[width=16cm]{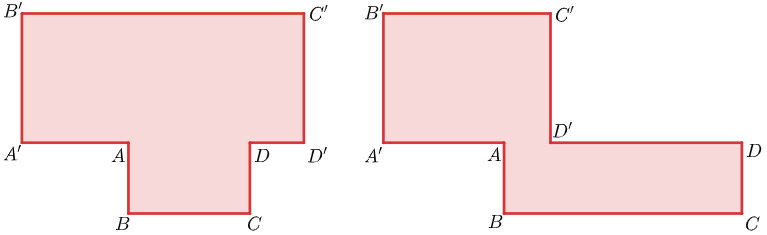}
	\caption{ T-shaped and Z-shaped octagons.}
	\label{TZ}
\end{figure}

Let $ABCD$ and {$A'B'C'D'$} be rectangles such that $A$ and $D$ lie on {$A'D'$} and the intersection of the rectangles is {$AD$}. Then the union of these rectangles is called \textit{a T-shaped octagon} (Figure~\ref{TZ}).

Let $ABCD$ and {$A'B'C'D'$} be rectangles such that $A$ and $D$ lie on {$A'D'$} and the intersection of the rectangles is {$AD'$} or {$A'D$}. Then the union of these rectangles is called \textit{a Z-shaped octagon} (Figure~\ref{TZ}).

Theorem~\ref{final_theorem} is true not only for $\Pi$-shaped octagons, but also for T-shaped and Z-shaped octagons. We prove it in this generality using Theorem~\ref{equivalence H or I} and the following lemma.

\begin{lemma} (Corollary from Lemma~\ref{correspondence-general}.) \label{correspondence-general2}  Let $P$ be a generic simple orthogonal octagon with $4$ vertical sides of signed lengths ${I_1, I_2, I_3, I_4}$ and ${x\mbox{-co}}$ordinates ${U_1, U_2, U_3, U_4}$. If $I_1+I_4=0$, then the following 2 conditions are
	equivalent:
	\begin{enumerate}[(1)]
		\item \label{1-lemma2} the octagon $P$ can be tiled by $m$ rectangles with ratios of the horizontal side to the vertical equal to $R_1, \dots, R_m$\textup{;}
		\item \label{2-lemma2} there is a planar two-port circuit with $m$ essential edges of resistances ${R_1,\dots,R_m}$, differences of incoming voltages $\Delta{U_{14}} := U_1-U_4$ and $\Delta{U_{23}} := U_2-U_3$, voltage drop value $U_1-U_2$, and incoming currents ${I_1, I_2}$.
	\end{enumerate}
In particular, all such two-port networks (corresponding to the two-port circuits) for a given octagon P are $\Pi$-equivalent.
\end{lemma}
\begin{proof}
$(\ref{1-lemma2})\Rightarrow(\ref{2-lemma2})$. 
Assume (\ref{1-lemma2}). Then by Lemma~\ref{correspondence-general2} there is a planar two-port circuit with $4$ terminals, $m$ essential edges of resistances ${R_1,\dots,R_m}$, incoming voltages ${U_1, U_2, U_3, U_4}$, and incoming currents ${I_1, I_2, I_3, I_4}$. By the assumptions of Lemma~\ref{correspondence-general2}, we have $I_1+I_4=0$. Then this is a two-port circuit and we get (\ref{2-lemma2}).

$(\ref{2-lemma2})\Rightarrow(\ref{1-lemma2})$. 
Consider the planar electrical circuit with incoming voltages ${U_1, U_2, U_3, U_4}$, with the same edges and the same resistances as the two-port circuit from (\ref{2-lemma2}). The currents in both circuits are the same. Then by Lemma~\ref{correspondence-general} we get (\ref{1-lemma2}).
\end{proof}

\begin{proof}[Proof of Theorem~\ref{final_theorem}]
	Fix an enumeration of the sides of the octagon, starting with $B'A'$. Then for a $\Pi$-, T- or Z-shaped octagon the equality $I_1+I_4=0$ is satisfied. By Lemma~\ref{correspondence-general2}, if such an octagon is tiled by $m$ squares, then we have (\ref{2-lemma2}).
	
By Theorem~\ref{equivalence H or I} and Remark~\ref{rational}, this planar two-port network (corresponding to the two-port circuit) is $\Pi$-equivalent to a network with no more than $5$ edges of rational resistances. By the $\Pi$-equivalence, it follows that the resulting circuit with the differences of incoming voltages $\Delta{U_{14}}$ and $\Delta{U_{23}}$ has the same voltage drop $U_1-U_2$ and the same incoming currents ${I_1 , I_2}$. Then by Lemma~\ref{correspondence-general2} the octagon can be tiled by no more than 5 rectangles with rational aspect ratios, because the resistances of the edges are rational.
\end{proof}

Under the assumptions of Theorem~\ref{final_theorem}, the number 5 is minimal.

\begin{figure}[h]
	\center
	\includegraphics[width=9cm]{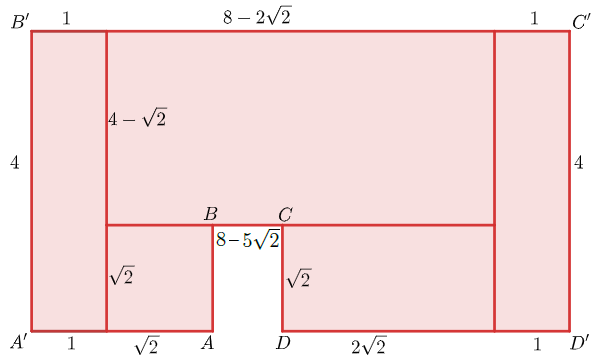}
	\caption{A $\Pi$-shaped octagon that cannot be tiled by less than 5 rectangles with rational aspect ratios.}
	\label{example2}
\end{figure}

\begin{example}\label{example_tiling>4} (Figure~\ref{example2})
The $\Pi$-shaped octagon with $A'B'=4, B'C'=10-2\sqrt{2}, AA'=1+\sqrt{2}, AB=\sqrt{2}, DD'=1+2\sqrt{2}$ can be tiled by 5 rectangles with rational aspect ratios, but cannot be tiled by less than 5 rectangles with rational aspect ratios.
\end{example}

\begin{proof} The tiling by 5 rectangles is shown in Figure~\ref{example2}. By Lemma~\ref{correspondence-general2}, there exists a planar two-port circuit with $5$ edges $15, 25, 36, 46,$ and $56$ of rational resistances such that $I_1 = 4, I_2 = -\sqrt{2}, I_3 = \sqrt{2}, I_4 = -4, U_1-U_2=-(1+\sqrt{2}),  U_3-U_4=-(1+2\sqrt{2}), U_1-U_4=-(1+\sqrt{2}),  U_3-U_4=-(1+2\sqrt{2}), \Delta{U_{14}} := U_1-U_4,$ and $\Delta{U_{23}} := U_2-U_3$. This circuit satisfies the assumptions of Proposition~\ref{contrexample}.

Assume the converse: this octagon can be tiled by less than 5 rectangles with rational aspect ratios. Then, similarly to the previous paragraph, by Lemma~\ref{correspondence-general2} we get a two-port circuit with less than 5 edges and the same differences of the incoming voltages, the same voltage drop, and the same response. This contradicts to Proposition~\ref{contrexample}.
\end{proof}
	
\section{Acknowledgments}
I am grateful to my scientific advisor Mikhail Skopenkov for valuable discussions and help in writing and translating the paper.
Also I would like to thank Pasha Pylyavskyy, who read this paper and provided ideas for further research.

\bibliographystyle{elsarticle-num}

\end{document}


\title{Преобразования двухпортовых цепей и разбиение на прямоугольники}
	\author{Светлана Широковских}
	\date{}
	\maketitle	
	\long\def\shorten#1\endshorten{}	
	\shorten	
\begin{center}
	{\scshape Федеральное государственное автономное\\
		образовательное учреждение высшего образования\\
		<<Национальный исследовательский университет\\
		<<Высшая школа экономики>>\\[1ex]
		Факультет математики\par}
	
	\par\vfill
	
	\textbf{\large Широковских Светлана Сергеевна}
	
	\vspace{1.5cm}
	
	{\Large\bfseries
		Преобразования двухпортовых цепей и разбиение на прямоугольники
		\par}
	
	\vspace{1.5cm}
	
	Курсовая работа студентки 2 курса\\[1ex]
	образовательной программы бакалавриата <<Математика>>
	\par\vfill
	\noindent\hspace{0.52\textwidth}\parbox[t]{0.48\textwidth}{%
		Научный руководитель:\\[3pt]
		кандидат физико-математических наук,\\
		доцент\\
		Скопенков Михаил Борисович\\[2ex]
	}%
	\par\vfill
	Москва 2021
\end{center}
\thispagestyle{empty}
\pagebreak
\endshorten

\begin{abstract}
В данной работе исследуются двухпортовые электрические цепи, вводятся новые понятия перепада и П-эквивалентности. Основной результат --- каждая плоская двухпортовая цепь П-эквивалентна цепи с не более чем 5 рёбрами. С его помощью доказывается, что если восьмиугольник в форме буквы П разрезается на квадраты, то он разрезается не более чем на 5 прямоугольников с рациональным отношением сторон. 
Теорема Кеньона 1998 года позволяла доказать это только для 6 прямоугольников. 
\end{abstract}

\keywords{Электрическая цепь, двухпортовая цепь, разрезание, прямоугольник, эквивалентные преобразования}

\section{Введение}
\begin{figure}[h]
	\center
	\includegraphics[width=6cm]{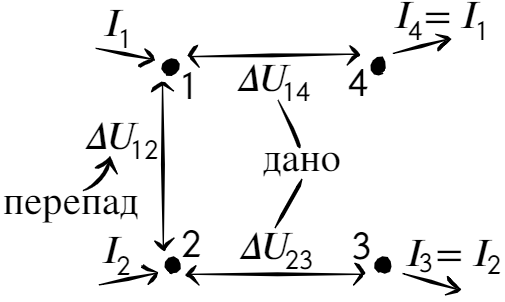}
	\caption{Двухпортовая цепь (см. раздел~\ref{prelim2}).}
	\label{fig1}
\end{figure}

В данной работе рассматриваются двухпортовые электрические цепи --- взвешенные графы с четырьмя выделенными вершинами, разбитыми на две пары, называющиеся \emph{портами}. Наглядным примером такой цепи может служить компьютер с двумя USB-портами или фонарик с двумя батарейками. В двухпортовых цепях задаются не потенциалы выделенных вершин, а разности потенциалов для каждого порта. При этом мы требуем, что если в одну вершину из порта втекает какой-то ток, то из второй вершины должен вытекать такой же (рисунок~\ref{fig1}; точное определение приводится в разделе~\ref{prelim2}).

В этой работе исследуется новое (хотя и естественное) понятие \emph{перепада напряжения} между портами, а также новое отношение \emph{П-эквивалентности}, которое требует не только одинаковых входящих токов в двух цепях, но и одинакового перепада напряжения между портами.

Эти понятия интересны с инженерной точки зрения. Если при преобразованиях цепи следить только за входящими токами, а не за сохранением перепада, то разность потенциалов между вершинами из разных портов может стать очень большой. Неформально, может произойти короткое замыкание, нарушающее нормальную работу цепи. Например, компьютер может сгореть, если перепаять его, следя только за сохранением напряжений и токов через USB-порты, но не перепада. П-эквивалентные преобразования помогают избежать таких проблем.

Основной результат работы состоит в том, что каждая плоская двухпортовая цепь П-эквивалентна цепи с не более чем 5 рёбрами (теорема~\ref{equivalence H or I}), причём число 5 минимально (пример~\ref{example_min=5}). Интересно, что если требовать только сохранения входящих токов, а не сохранения перепада напряжения, то число рёбер можно уменьшить до 3 (замечание~\ref{equivalence 3}).

\begin{figure}[h!]
	\includegraphics[width=15cm]{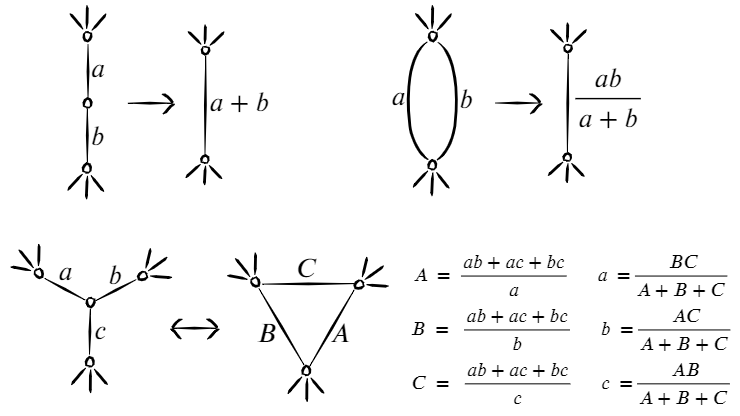}	
	\caption{Элементарные преобразования (см. раздел~\ref{prelim1}).}
	\label{elementary}
\end{figure}

\begin{figure}[h]
	\center
	\includegraphics[width=11cm]{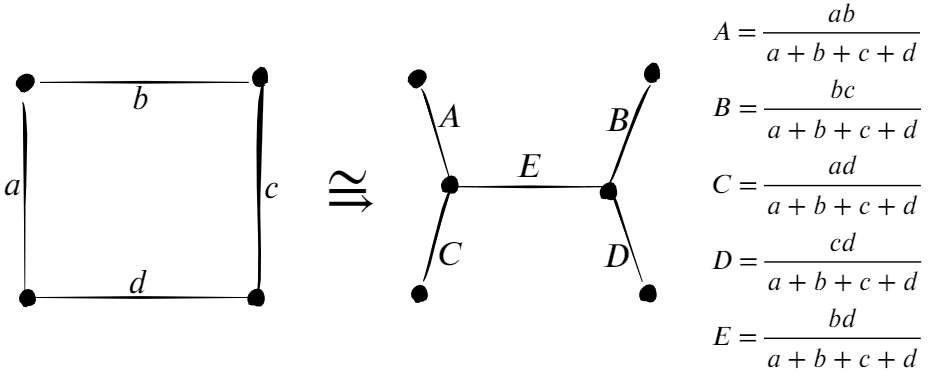}
	\caption{Box-H преобразование (см. следствие~\ref{Box-H transformation}).}
	\label{Box-H_transformation}
\end{figure}

Доказательство основано на том, чтобы при помощи П-эквивалентных преобразований уменьшить число ребер в цепи. Самые простые из таких преобразований хорошо известны; они изображены на рисунке~\ref{elementary} (см.~\cite{K},~\cite{CM},~\cite{ZSU}). В этой работе вводится новое \textit{Box-H преобразование} (рисунок~\ref{Box-H_transformation}).

Новые результаты о перепаде напряжения и П-эквивалентности позволяют решать новые задачи о замощениях (благодаря известной физической интерпретации замощений Р.~Л.~Брукса, К.~А.~Б.~Смита, А.~Г.~Стоуна и У.~Т.~Татта~\cite{BSST}, подробно описанной в~\cite{PS}). С их помощью мы доказываем, что если восьмиугольник в форме буквы П можно разрезать на квадраты, то его можно разрезать не более чем на $5$ прямоугольников с рациональным отношением сторон (теорема~\ref{final_theorem}). Известная ранее теорема Кеньона (теорема~\ref{Kenyon}) позволяла доказать это только для $6$ прямоугольников.

\subsection{Преобразования двухпортовых цепей}

Преобразования двухпортовых цепей, сохраняющие входящие токи, изучались в литературе (\cite{BIK},~\cite{Bess},~\cite{B}). Приведём пару примеров.

\begin{figure}[h]
	\center
	\includegraphics[width=6cm]{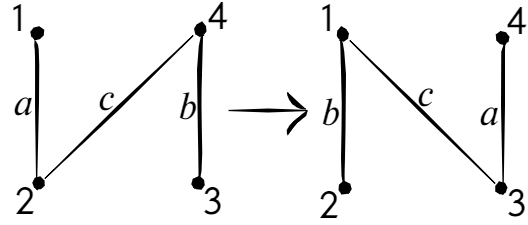}
	\caption{Двухпортовые цепи с одинаковыми входящими токами .}
	\label{N-N}
\end{figure}

Например, легко видеть, что при одинаковых разностях входных потенциалов двухпортовые цепи, изображённые на рисунке~\ref{N-N}, будут иметь одинаковые входящие токи. Действительно, если в правой цепи поменять местами вершины в каждом порту ($\{1, 4\}$ и $\{2, 3\}$), то получится левая цепь, однако разности входных потенциалов поменяют знак. Поэтому ток, выходящий из вершины 1 левой цепи, равен току, входящему в вершину 4 правой цепи. Рассмотрим порт $\{1, 4\}$ в левой цепи. В двухпортовых цепях мы требуем, что если в одну вершину пары втекает какой-то ток, то из второй вершины должен вытекать такой же. Поэтому в левой цепи ток, выходящий из вершины 1, равен току, входящему в вершину 4. По транзитивности ток, входящий в вершину 4, в левой и правой цепях одинаковый. Аналогично с остальными вершинам.

Будем называть двухпортовую цепь \textit{симметричной}, если существует автоморфизм взвешенного графа, переставляющий вершины в каждой паре $\{1, 2\}$ и $\{3, 4\}$.

\begin{figure}[h]
	\center
	\includegraphics[width=11cm]{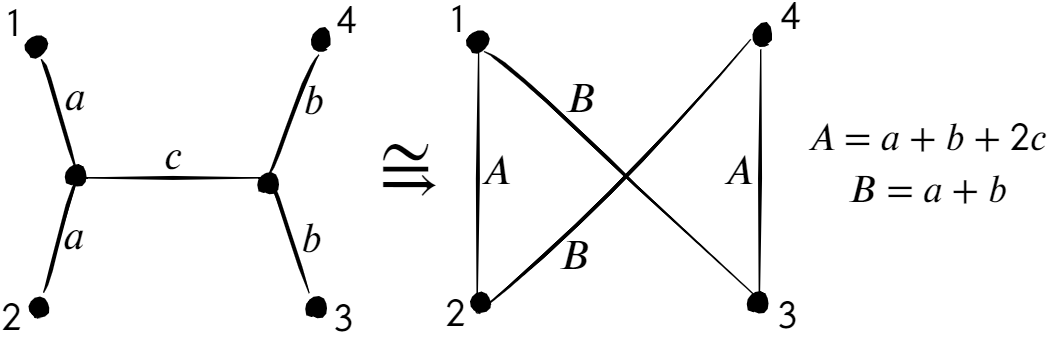}
	\caption{Пример двухпортовых цепей с одинаковыми входящими токами.}
	\label{wiki}
\end{figure}

В 1927 году Барлетт~\cite{B} доказал, что для любой симметричной двухпортовой цепи существует двухпортовая цепь с четырьмя ребрами ${12}, {34}, {13}$ и ${24}$ с такими же входящими токами (рис.~\ref{wiki}).

\subsection{Разрезания}

Приведём небольшой обзор известных результатов о разрезаниях на квадраты.

Естественный вопрос: какие многоугольники можно разрезать на квадраты? Ответ для прямоугольника был получен Максом Деном в 1903 году.

\begin{theorem} [Ден, 1903\label{dehn}] \textup{\cite{D}} Прямоугольник разрезается на квадраты \textup{(}не обязательно равные\textup{)} тогда и только тогда, когда отношение его сторон рационально.
\end{theorem}

Здесь и далее \textit{разрезанием многоугольника на прямоугольники} называется представление этого многоугольника в виде объединения конечного числа замкнутых прямоугольников, внутренности которых не пересекаются.

Следующий результат в этом направлении был получен Ричардом Кеньоном спустя почти сто лет.

\begin{theorem}[Кеньон, 1998] \label{th-kenyon}
	\textup {(Ср.~\cite[Theorem~9]{K})} Пусть $A_1A_2A_3A_4A_5A_6$ — шестиугольник с прямыми углами, вершины которого пронумерованы против часовой стрелки, начиная с вершины невыпуклого угла. Шестиугольник $A_1A_2A_3A_4A_5A_6$ можно разрезать на квадраты, если и только если система линейных уравнений
	\begin{equation}
		\label{syst}
		\begin{cases}
			A_3A_4\cdot x+ A_1A_2\cdot y&=A_2A_3,\\
			A_5A_6\cdot z-A_1A_2\cdot y&=A_6A_1;
		\end{cases}
	\end{equation}
	имеет решение в неотрицательных рациональных числах $x,y,z$.
\end{theorem}

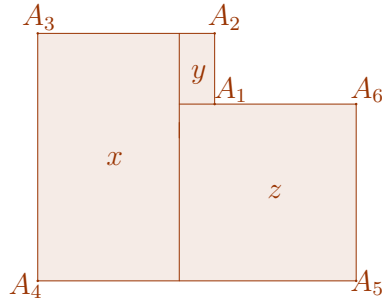
\begin{figure}[htbp]
	\center
	\definecolor{zzttqq}{rgb}{0.6,0.2,0}
\begin{tikzpicture}[line cap=round,line join=round,>=triangle 45,x=0.5cm,y=0.5cm]
\clip(-4.88,-1.69) rectangle (5.72,6.85);
\fill[color=zzttqq,fill=zzttqq,fill opacity=0.1] (1,6) -- (-4,6) -- (-4,-1) -- (5,-1) -- (5,4) -- (1,4) -- cycle;
\draw [color=zzttqq] (1,6)-- (-4,6);
\draw [color=zzttqq] (-4,6)-- (-4,-1);
\draw [color=zzttqq] (-4,-1)-- (5,-1);
\draw [color=zzttqq] (5,-1)-- (5,4);
\draw [color=zzttqq] (5,4)-- (1,4);
\draw [color=zzttqq] (1,4)-- (1,6);
\draw [color=zzttqq] (0,6)-- (0,3.04);
\draw [color=zzttqq] (0,-1)-- (0,3.48);
\draw [color=zzttqq] (0,4)-- (1,4);
\draw[color=zzttqq] (-2.35,2.94) node[anchor=north west] {$x$};
\draw[color=zzttqq] (0.05,5.4) node[anchor=north west] {$y$};
\draw[color=zzttqq] (2.2,2) node[anchor=north west] {$z$};
\fill [color=zzttqq] (1,6) circle (0.5pt);
\draw[color=zzttqq] (1.28,6.45) node {$A_2$};
\fill [color=zzttqq] (-4,6) circle (0.5pt);
\draw[color=zzttqq] (-3.96,6.45) node {$A_3$};
\fill [color=zzttqq] (-4,-1) circle (0.5pt);
\draw[color=zzttqq] (-4.34,-1.14) node {$A_4$};
\fill [color=zzttqq] (5,-1) circle (0.5pt);
\draw[color=zzttqq] (5.32,-1.06) node {$A_5$};
\fill [color=zzttqq] (5,4) circle (0.5pt);
\draw[color=zzttqq] (5.32,4.4) node {$A_6$};
\fill [color=zzttqq] (1,4) circle (0.5pt);
\draw[color=zzttqq] (1.46,4.4) node {$A_1$};
\end{tikzpicture} 
	\caption{\textup{\cite{PS}} Шестиугольник $A_1A_2A_3A_4A_5A_6$.}
	\label{ell}
\end{figure}

Рисунок~\ref{ell} показывает, почему выполняется часть <<если>>. Если система линейных уравнений~\eqref{syst} имеет решение в неотрицательных рациональных числах $x,y,z$, то шестиугольник можно разрезать на 3 прямоугольника (возможно, вырожденных) с рациональными отношениями сторон $x,y,z$, как показано на рисунке. А каждый такой прямоугольник можно разрезать на квадраты.

Многоугольник называется \emph{ортогональным}, если его стороны параллельны координатным осям. Будем говорить, что он \emph{общего положения}, если координаты вертикальных сторон по горизонтальной оси попарно различны.

\begin{theorem} \label{corollary_about_6}
	(следствие из теоремы~\ref{Kenyon} и леммы~\ref{correspondence-general} ниже; Кеньон, 1998) Пусть P — ортогональный $2n$-угольник общего положения. Если его можно разрезать на квадраты, то его можно разрезать не более чем на $n(n-1)/2$ прямоугольников с рациональными отношениями сторон. В частности, ортогональный шестиугольник можно разрезать на 3 прямоугольника с рациональными отношениями сторон, а восьмиугольник -- на 6.
\end{theorem}

В данной работе для П-образного восьмиугольника доказывается следующий более сильный факт.

\begin{figure}[h]
	\center
	\includegraphics[width=10cm]{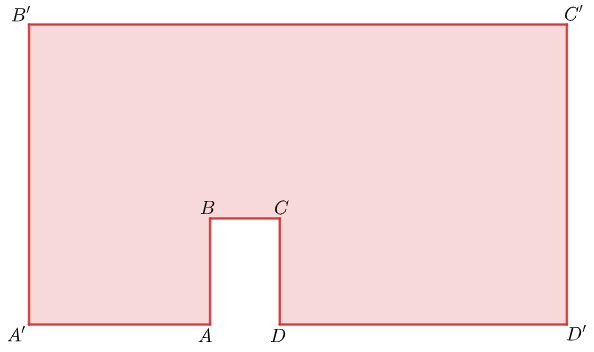}
	\caption{ П-образный восьмиугольник.}
	\label{P-8}
\end{figure}

Пусть есть 2 прямоугольника $ABCD$ и {$A'B'C'D'$}, причем $A$ и $D$ лежат на отрезке {$A'D'$}, а $B$ и $C$ лежат строго внутри {$A'B'C'D'$}. Тогда замыкание дополнения {$A'B'C'D'\setminus ABCD$} называется \textit{П-образным восьмиугольником} (рис.~\ref{P-8}).

\begin{theorem}\label{final_theorem}
	Если П-образный восьмиугольник можно разрезать на квадраты, то его можно разрезать на 5 или меньше прямоугольников с рациональными отношениями сторон.
\end{theorem}

\subsection{Структура работы}
В разделе~\ref{prelim} приводятся известные определения и теоремы об электрических цепях. В подразделе~\ref{prelim1} вводится строгое определение электрической цепи, а также определения ответа и эквивалентности. Кроме того вводятся некоторые эквивалентные преобразования электрических цепей. В подразделе~\ref{prelim2} вводится определение двухпортовой электрической цепи, а понятие ответа обобщается на случай двухпортовых цепей. В коротком подразделе~\ref{prelim3} приводится лемма, связывающая разрезания и электрические цепи.

В разделе~\ref{new_equivalence} вводятся новые понятия перепада и П-эквивалентности. Затем доказывается П-эквивалентность Box-цепи и H-цепи (рис.~\ref{Box-H_transformation}). Далее следует основная теорема данной работы (теорема~\ref{equivalence H or I}) и её доказательство. После этого следует пример, доказывающий, что нельзя усилить формулировку теоремы.

В разделе~\ref{tiling_octagons} при помощи теоремы~\ref{equivalence H or I} доказывается теорема~\ref{final_theorem} о разрезании П-образного восьмиугольника и приводится пример такого восьмиугольника, который можно разрезать минимум на 5 прямоугольников с рациональным отношением сторон.

\section{Предварительные сведения}\label{prelim}
\subsection{Электрические цепи}\label{prelim1}
\textit{Электрической цепью с ${t}$ входами} называется связный граф, каждому ребру которого сопоставлено некоторое положительное действительное число (\textit{проводимость}) и ${t}$ вершин которого отмечены (\textit{граничные вершины}). Для простоты будем считать, что граф не содержит ни кратных рёбер, ни петель. Обобщения для графов с кратными рёбрами оставляем читателю.

Электрическую цепь будем называть \textit{плоской}, если она расположена в единичном круге так, что граничные вершины цепи принадлежат границе круга и никакие два ребра не имеют общих точек, кроме, возможно, их общей вершины.

Зафиксируем нумерацию вершин цепи числами ${1}$, ${2}$, ${\dots}$, ${n}$ так, чтобы вершины ${1, \dots, t}$ были бы граничными. В случае если цепь плоская, будем нумеровать граничные вершины против часовой стрелки по границе круга. Обозначим ребро между вершинами $k$ и $l$ через ${kl}$, проводимость ребра ${kl}$ --- через ${c_{kl}}$. Положим ${c_{kl}:=0}$, если между вершинами ${k}$ и ${l}$ нет ребра.

\textit{Замкнутой электрической цепью} называется электрическая цепь вместе с ${t}$ действительными числами ${U_1, \dots, U_t}$ (\textit{входными потенциалами}), сопоставленными граничным вершинам цепи.

Каждой вершине ${k}$ (${1\le k\le n}$) замкнутой электрической цепи сопоставим число ${U_k}$ --- \textit{потенциал} вершины, а каждой (упорядоченной) паре вершин ${1\le k,l\le n}$ --- число ${I_{kl}}$, называемое \textit{током} через ребро ${kl}$ в направлении от $k$ к $l$ (если между $k$ и $l$ нет ребра, то это число будет нулём по аксиоме ${(C)}$ ниже).
Эти числа определяются следующими двумя аксиомами:

\begin{enumerate}[(1)]
\item[(C)]\textit{Закон Ома.} Для каждой пары вершин ${k,l}$ выполнено равенство ${I_{kl}=c_{kl}(U_k-U_l)}$.
\item[(I)]\textit{Правило Кирхгофа.} Для каждой вершины ${k>t}$ выполнено равенство ${{\sum_{l=1}^{n} I_{kl}=0}}$.
\end{enumerate}

Числа ${U_k}$ и ${I_{kl}}$ однозначно определяются этими аксиомами, что утверждает следующая теорема.

\begin{theorem} \textup{~\cite[Теорема~2.1.]{PS}} 
Для любой замкнутой электрической цепи
система линейных уравнений \textup{(C),\7(I)} от переменных ${U_k}$, где ${t< k\le n}$, и ${I_{kl}}$, где ${1\le k,l\le n}$, имеет единственное решение.
\end{theorem}

Величина, обратная к проводимости ребра, называется его \emph{сопротивлением} $R_{kl} := 1/c_{kl}$.
 Числа $(I_{1}, \dots, I_{t}) := ({{\sum_{k=1}^{n} I_{1k}}}, \dots, {{\sum_{k=1}^{n} I_{tk}}})$ называются \textit{входными токами}. \textit{Ответом} электрической цепи называется отображение $\mathds{R}^t \to \mathds{R}^t$: $({U_{1}}, \dots, {U_{t}}) \mapsto (I_{1}, \dots, I_{t})$. Очень интересно описание ответов плоских цепей (см.~\cite{CM},~\cite{K}). Две цепи \textit{эквивалентны}, если их ответы равны.

Результаты этого раздела хорошо известны и легко выводятся из определений (\cite[Теорема об электрических преобразованиях, \S2.3]{SSU}).

\begin{figure}[h]
	\center
	\includegraphics[width=8cm]{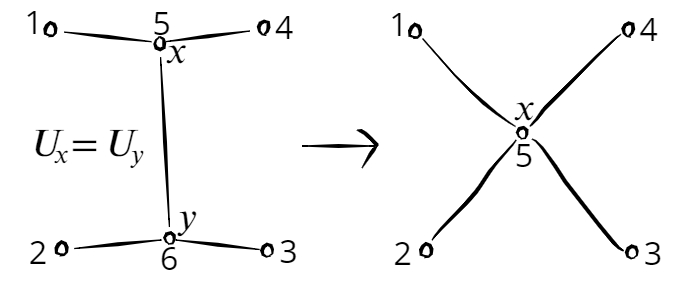}
	\caption{Объединение вершин $x$ и $y$ с одинаковым потенциалом в I-цепи.}
	\label{same_potentials2}
\end{figure}

\begin{predl}\label{same potential0} (рис.~\ref{same_potentials2})
Пусть $xy$ --- такое ребро электрической цепи, что $U_x=U_y$ при любых входных потенциалах. Тогда если заменить ребра $zy$ для всех $z \neq x$ на ребра $zx$ с такими же сопротивлениями, а затем удалить ребро $xy$ и вершину $y$, то получится эквивалентная цепь.
\end{predl}

Такое преобразование цепи назовём \emph{объединением вершин $x$ и $y$ с одинаковым потенциалом}.

\begin{remark}\label{remark}(Параллельное соединение, рис.~\ref{elementary}, справа вверху).
Если у цепи есть два кратных ребра сопротивлениями $a$ и $b$, то можно заменить их на одно ребро сопротивлением $ab/(a+b)$ между теми же вершинами. Получившаяся цепь будет эквивалентна изначальной. При помощи таких параллельных соединений из любой цепи можно получить эквивалентную ей цепь без кратных рёбер, причём число рёбер в цепи уменьшится. И если сопротивления рёбер в изначальной цепи были рациональными, то и в новой цепи сопротивления будут рациональны. Это позволяет нам ограничиться рассмотрением цепей без кратных рёбер.
\end{remark}

\begin{predl}\label{series equivalence}(рис.~\ref{elementary}, слева вверху)
Для любых $a, b > 0$ следующие две цепи с двумя входами эквивалентны:
\begin{enumerate}[(1)]
\item[(1)] Цепь на трех вершинах с двумя ребрами 13 и 23 сопротивлениями 
\begin{align*}
R_{13} = a, R_{23} = b.
\end{align*} 
\item[(2)] Цепь на двух вершинах с одним ребром 12 сопротивлением
\begin{align*}
R_{12} = a+b.
\end{align*}
\end{enumerate}
\end{predl}

\begin{corollary}\label{series transformation} (рис.~\ref{elementary}, слева вверху)
Пусть в электрической цепи есть неграничная вершина $y$, из которой выходят ровно два ребра $xy$ и $yz$. Тогда если заменить эти два ребра на одно ребро $xz$ сопротивлением $R_{xz}:=R_{xy}+R_{yz}$, а вершину $y$ удалить, то получится эквивалентная цепь.
\end{corollary}
	
Такое преобразование цепи назовём \emph{последовательным соединением} рёбер $xy$ и $yz$.

\begin{predl}\label{Delta-Y equivalence} (рис.~\ref{elementary}, снизу)
Для любых $a, b, c > 0$ следующие две цепи с тремя входами эквивалентны:
\begin{enumerate}[(1)]
\item[($\Delta$)] Цепь на трех вершинах с тремя ребрами 12, 23, 13 сопротивлениями
\begin{align*}
R_{12} = a, R_{23} = b, R_{13} = c.
\end{align*}
		
\item[($Y$)] Цепь на четырех вершинах с тремя рёбрами 14, 24, 34 сопротивлениями
\begin{align*}
R_{14} = ac/(a+b+c), R_{24} = ab/(a+b+c), R_{34} = bc/(a+b+c).
\end{align*}
\end{enumerate}
\end{predl}

\begin{corollary}\label{Delta-Y transformation} (рис.~\ref{elementary}, снизу)
Пусть в электрической цепи есть три попарно соединенных ребрами вершины $x$, $y$ и $z$. Тогда если добавить в цепь новую неграничную вершину $w$ и заменить ребра $xy$, $xz$, $yz$ на ребра $xw, yw, zw$ сопротивлениями 
\begin{align*} &R_{xw}:=R_{xy}R_{xz}/(R_{xy}+R_{xz}+R_{yz}),\\
&R_{yw}:=R_{xy}R_{yz}/(R_{xy}+R_{xz}+R_{yz}),\\
&R_{zw}:=R_{xz}R_{yz}/(R_{xy}+R_{xz}+R_{yz}),
\end{align*}
то получится эквивалентная цепь.
\end{corollary}

Такое преобразование цепи назовём \emph{преобразованием треугольник-звезда}, примененным к рёбрам $xy, yz, xz$. Обратное преобразование назовём преобразованием \emph{звезда-треугольник}. Преобразования параллельное соединение, последовательное соединение, преобразование треугольник-звезда и звезда-треугольник называются \emph{элементарными преобразованиями} (рис.~\ref{elementary}).

\begin{theorem}(Кеньон, 1998) \label{Kenyon} 
	\textup {~\cite[Theorem~8]{K}}
Для любой плоской электрической цепи $G$ с $t$ входами существует эквивалентная ей плоская цепь, которая имеет не более $t(t-1)/2$ ребер и получается из цепи $G$ последовательностью элементарных преобразований.	
\end{theorem}

\subsection{Двухпортовысли две электрические цепи}\label{prelim2}

\textit{Двухпортовой цепью} или \textit{четырёхполюсником} называется электрическая цепь с 4 входами, разбитыми на 2 пары, называющиеся \emph{портами}

Зафиксируем нумерацию вершин цепи числами ${1}$, ${2}$, ${\dots}$, ${n}$ так, чтобы вершины ${1}$, ${2}$, ${3}$, ${4}$ были бы граничными и разбитыми на пары $\{{1}, {4}\}$ и $\{{2}, {3}\}$. Как и раньше, обозначим проводимость ребра между вершинами ${k}$ и ${l}$ через ${c_{kl}}$. Положим ${c_{kl}:=0}$, если между вершинами ${k}$ и ${l}$ нет ребра. Плоская двухпортовая цепь определяется аналогично плоской электрической цепи. Если двухпортовая цепь плоская, то мы считаем, что граничные вершины пронумерованы против часовой стрелки вдоль границы круга.

\textit{Замкнутой двухпортовой цепью} называется двухпортовая цепь вместе с 2 действительными числами $\Delta{U_{14}}$ и $\Delta{U_{23}}$ (\textit{разности входных потенциалов}), сопоставленными парам граничных вершин цепи.

Каждой вершине ${k}$ (${1\le k\le n}$) двухпортовой цепи сопоставим число ${U_k}$ --- \textit{потенциал} вершины, а каждой (упорядоченной) паре вершин ${1\le k,l\le n}$ --- число ${I_{kl}}$, называемое \textit{током} через ребро $kl$ в направлении от $k$ к $l$ (${I_{kl}} = 0$, если между вершинами $k$ и $l$ нет ребра). Потенциалы и токи определяются следующими аксиомами:

\begin{enumerate}[(1)]
\item[(C)]\textit{Закон Ома.} Для каждой пары вершин ${k,l}$ выполнено равенство ${I_{kl}=c_{kl}(U_k-U_l)}$.
\item[(I)] \textit{Правило Кирхгофа.} Для каждой вершины ${k>4}$ выполнено равенство ${{\sum_{l=1}^{n} I_{kl}=0}}$.
\item[(P)] \textit{Условие изолированности портов.} ${{\sum_{l=1}^{n} I_{1l}+\sum_{l=1}^{n} I_{4l}=0}}$.
\item[(B)]\textit{Граничные условия.} $ U_1 - U_4 = \Delta{U_{14}}$ и $U_2 - U_3 = \Delta{U_{23}}$.
\end{enumerate}

\begin{theorem}\label{exist_uniq}
Для любой замкнутой двухпортовой цепи существуют единственные ${I_{kl}}$, ${1\le k,l\le n}$, и единственные с точностью до добавления общей постоянной числа ${U_{k}}$, ${1\le k\le n}$, удовлетворяющие свойствам (C), (I), (P), (B).
\end{theorem}

Доказательство теоремы~\ref{exist_uniq} хорошо известно, но мы приведём его для полноты. Оно аналогично доказательству теоремы~2.1 из \textup{~\cite{PS}} и использует следующую лемму:

\begin{lemma} \textup{~\cite[Лемма~5.1]{PS}} \label{energy} \label{fullpower} 
	Рассмотрим электрическую цепь с вершинами ${1, \dots, n}$, такую, что вершины ${1, \dots, t}$ --- граничные.
	Предположим, что числа ${U_k}$, где ${1\le k\le n}$, и ${I_{kl}}$, где ${1\le k,l\le n}$, удовлетворяют закону Ома~\textup{(C)} и правилу Кирхгофа~\textup{(I)} из \S\ref{prelim}. Положим ${I_u=\sum_{k=1}^n I_{uk}}$. Тогда
	$$\sum\limits_{1\le k<l\le n}(U_k-U_l)I_{kl}
	= \sum\limits_{1\le u\le t}U_u I_u.$$
\end{lemma}

\begin{proof}[Доказательство теоремы~\ref{exist_uniq}.] 
Легко видеть, что набор токов ${I_{kl}}$ и потенциалов ${U_{k}}$ удовлетворяет свойствам (C), (I), (P), (B) тогда и только тогда, когда им удовлетворяет набор токов ${I_{kl}}$ и потенциалов ${U_{k}} + C$. Будем доказывать существование и единственность решения с дополнительным условием $U_1=0$.

\textit{Единственность}.
Предположим, что два набора токов ${I^{I,II}_{kl}}$ и потенциалов ${U^{I,II}_k}$ подчиняются (C), (I), (P), (B) и ${U^{I}_1 = U^{II}_2 = 0}$ при некоторых разностях входных потенциалов ($\Delta{U_{14}}^{I}$ = $\Delta{U_{14}}^{II}$ и $\Delta{U_{23}}^{I}$ = $\Delta{U_{23}}^{II}$). Тогда их разность ${I_{kl}=I^{I}_{kl}-I^{II}_{kl}}$ и ${U_{k}=U^{I}_{k}-U^{II}_{k}}$ подчиняется (C), (I), (P), (B) при нулевых разностях входных потенциалов $\Delta{U_{14}}=0$ и $\Delta{U_{23}}=0$.
 Тогда $U_1 = U_4$ и $U_2 = U_3$. Из $(I)$ и $(P)$ следует, что $I_1=-I_4$ и $I_2=-I_3$. Получаем, что $\sum_{u=1}^4 U_u I_u = 0$.

По лемме~\ref{energy}: $$\sum\limits_{1\le k<l\le n}(U_k-U_l)^2 c_{kl}
=\sum\limits_{1\le k<l\le n}(U_k-U_l)I_{kl}
= \sum\limits_{1\le u\le 4}U_u I_u = 0.$$

Для всех ${k,l}$ либо ${c_{kl}>0}$, либо ${c_{kl}=0}$. Таким образом, каждое слагаемое $ c_{kl}(U_k-U_l)^2=0$. Поскольку граф связный, все потенциалы ${U_{k}}$ равны между собой. Но $U_1 = U_1^{I} - U_1^{II}=0$. Следовательно, ${U_k=0}$, ${I_{kl}=0}$, и поэтому ${I^{I}_{kl}=I^{II}_{kl}}$, ${U^{I}_{k}=U^{II}_{k}}$ для всех ${k,l}$.

\textit{Существование}. Количество уравнений в системе (C), (I), (P), (B) равно количеству неизвестных (при фиксированном $U_1=0$).
Мы доказали, что эта система имеет единственное решение при $\Delta{U_{14}}=0$, $\Delta{U_{23}}=0$. По конечномерной альтернативе Фредгольма она имеет решение для любых $\Delta{U_{14}}$, $\Delta{U_{23}}$.
\end{proof}

Числа $(I_{1}, I_{2}) := ({{\sum_{k=1}^{n} I_{1k}}}, {{\sum_{k=1}^{n} I_{2k}}})$ называются \emph{входными токами}. \textit{Ответом} двухпортовой цепи (или \textit{матрицей Y-параметров}~\cite{BIK}) называется отображение $\mathds{R}^2 \to \mathds{R}^2$: $(\Delta{U_{14}}, \Delta{U_{23}}) \mapsto (I_{1}, I_{2})$.

\begin{example}\label{Box_response}
	Двухпортовая цепь на четырёх вершинах с четырьмя рёбрами $12, 23, 34, 14$ сопротивлениями
	$R_{12} = a, R_{14} = b, R_{34} = c, R_{23} = d$ (\textit{Box-цепь}, см. левую цепь на рисунке~\ref{Box-H_transformation}) имеет следующую матрицу ответа: 
	\begin{center}
		${\left( \begin{matrix}
				(a+b+c)/(ab+ac) & -1/(a+c)\\
				-1/(a+c) & (a+c+d)/(ad+cd)
			\end{matrix}\right)}$.
	\end{center}
\end{example}

\begin{proof}
	Согласно аксиомам $(C), (I), (P), (B),$ система уравнений для цепи:
	\begin{equation}
		\label{syst2}
		\begin{cases}
		U_1-U_4=\Delta U_{14}, \\
		U_2-U_3=\Delta U_{23},\\
		I_{14}b=U_1-U_4,\\
		I_{12}a=U_1-U_2,\\
		I_{43}c=U_4-U_3,\\
		I_{23}d=U_2-U_3,\\ 
		I_{12}+I_{14}+I_{23}+I_{43}=0.
		\end{cases}
	\end{equation}
	Откуда получаем:
	\begin{align*}
		&I_{14}=\Delta{U_{14}}/b,\\
		&I_{12}=(\Delta{U_{14}}-\Delta{U_{23}})/(a+c),\\
		&I_{43}=(\Delta{U_{23}}-\Delta{U_{14}})/(a+c),\\
		&I_{23}=\Delta{U_{23}}/d,\\
		&I_{1} = I_{14}+I_{12}=((a+b+c)\Delta{U_{14}}-b\Delta{U_{23}})/(b(a+c)),\\
		&I_{2} = I_{21}+I_{23}=((a+c+d)\Delta{U_{23}}-d\Delta{U_{14}})/(d(a+c)).
	\end{align*}
	Следовательно:
	\begin{center}
		${\left( \begin{matrix}
				(a+b+c)/(ab+ac) & -1/(a+c)\\
				-1/(a+c) & (a+c+d)/(ad+cd)
			\end{matrix}\right)} \cdot
		{\left( \begin{matrix}
				\Delta{U_{14}}\\
				\Delta{U_{23}}
			\end{matrix}\right)} = 
		{\left( \begin{matrix}
				{I_{1}}\\
				{I_{2}}
			\end{matrix}\right)}
		$.
		
	\end{center}
	
\end{proof}

\subsection{Электрические цепи и разрезания}\label{prelim3}

\begin{figure}[h]
	\center
	\includegraphics[width=12cm]{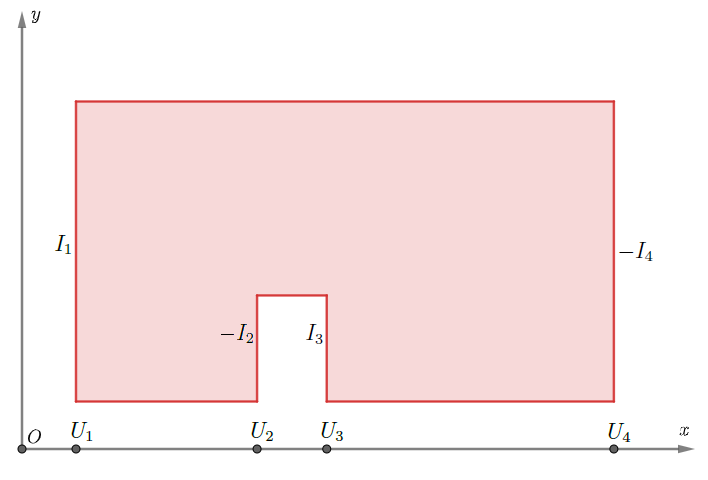}
	\caption{Ортогональный многоугольник.}
	\label{example_ort}
\end{figure}

Пусть $P$ --- \emph{ортогональный} многоугольник (рис.~\ref{example_ort}), то есть многоугольник со сторонами, параллельными координатным осям. Будем считать, что многоугольник $P$ --- \emph{простой}, то есть его граница $\partial P$ состоит из одной связной компоненты. Занумеруем стороны многоугольника, параллельные оси $y$, против часовой стрелки в порядке их следования на границе $\partial P$.
Обозначим через $t$ число этих сторон.
Обозначим через $I_u$ \emph{ориентированную длину} неориентированной стороны $u$, то есть длину стороны $u$ со знаком \2$+$\3 (соответственно, \2$-$\3), если многоугольник $P$ локально лежит правее (соответственно, левее) неориентированной стороны $u$. Обозначим через $U_u$ координату стороны $u$ по оси $x$. Будем говорить, что $P$ --- \emph{общего положения}, если числа ${U_1, \dots, U_t}$ попарно различны.

 Мы будем говорить, что ребро $kl$ замкнутой цепи \emph{существенно}, если $I_{kl} \neq 0$.

\begin{lemma}\textup{~\cite[Лемма~5.4.]{PS}} \label{correspondence-general} Пусть $P$ --- ортогональный многоугольник общего положения с $t$ вертикальными сторонами, имеющими ориентированные длины ${I_1, \dots, I_t}$ и ${x\mbox{-ко}}$орди\-на\-ты ${U_1, \dots, U_t}$. Тогда следующие $2$ условия эквивалентны\textup{:}
	\begin{enumerate}[(1)]
		\item \label{c-g-1} многоугольник $P$ можно разрезать на $m$ прямоугольников с отношениями горизонтальной стороны к вертикальной, равными $R_1, \dots, R_m$\textup{;}
		\item \label{c-g-2} существует плоская замкнутая электрическая цепь с $t$ граничными вершинами, состоящая из $m$ существенных рёбер сопротивлениями ${R_1,\dots,R_m}$, с входными потенциалами ${U_1, \dots, U_t}$ и входными токами ${I_1, \dots, I_t}$.
	\end{enumerate}
\end{lemma}

\section{П-эквивалентность двухпортовых электрических цепей}\label{new_equivalence}
Теперь введём новое понятие, которое появляется в работе впервые. \textit{Перепадом} двухпортовой цепи назовём отображение $\mathds{R}^2 \to \mathds{R}$: $(\Delta{U_{14}}, \Delta{U_{23}}) \mapsto U_{1}-U_{2}$. Будем называть двухпортовые цепи \textit{П-эквивалентными}, если их ответы и перепады равны. В лемме~\ref{correspondence-general2} мы увидим, как эти понятия естественно возникают при изучении разрезаний на прямоугольники.

\begin{example}\label{Box_drop}
У Box-цепи (см. пример~\ref{Box_response} и рис.~\ref{Box-H_transformation}) перепад $U_{1}-U_{2} = (\Delta{U_{14}}-\Delta{U_{23}})a/(a+c)$.
\end{example}

\begin{proof}
Это следует непосредственно из~(\ref{syst2}).
\end{proof}

Отметим следующие очевидные утверждения.

\begin{predl}\label{same potential} (Рис.~\ref{same_potentials2}, ср. предложение~\ref{same potential0})
	Пусть $xy$ --- такое ребро двухпортовой цепи, что $U_x=U_y$ при любых разностях входных потенциалов. Тогда если заменить ребра $zy$ для всех $z \neq x$ на ребра $zx$ с такими же сопротивлениями, а затем удалить ребро $xy$ и вершину $y$, то получится П-эквивалентная цепь.
\end{predl}

Данное преобразование, обобщающее объединение вершин с одинаковым потенциалом, по-прежнему называется \emph{объединением вершин $x$ и $y$ с одинаковым потенциалом}.

\begin{predl}\label{equel=>}
Если две электрические цепи c четырьмя входами эквивалентны, то двухпортовые цепи с такими же рёбрами с такими же сопротивлениями будут П-эквивалентны.
В частности, элементарные преобразования сохраняют классы П-эквивалентности двухпортовых цепей.
\end{predl}

Нашим основным инструментом будет новое преобразование двухпортовых цепей, сохраняющее как ответ, так и перепад (рис.~\ref{Box-H_transformation}). Это преобразование, в частности, показывает, что утверждение, обратное к предложению~\ref{equel=>} неверно (можно проверить непосредственно, вычислив матрицы ответов, что это не эквивалентное преобразование).
\textit{H-цепью} назовём двухпортовую цепь на 6 вершинах с 5 рёбрами: ${15}, {25}, {46}, {36}, {56}$.

\begin{theorem}\label{Box-H}(рис.~\ref{Box-H_transformation})
Для любых $a, b, c, d > 0$ следующие две двухпортовые цепи П-эквивалентны:
\begin{enumerate}[(1)]
\item[(Box)] Цепь на четырёх вершинах с четырьмя рёбрами $12, 23, 34, 14$ сопротивлениями
\begin{align}
\nonumber
R_{12} = a, R_{14} = b, R_{34} = c, R_{23} = d.	\end{align}
		
\item[(H)] Цепь на шести вершинах с пятью рёбрами $15, 25, 36, 46, 56$ сопротивлениями
\begin{align*}
& R_{15} = ab/(a+b+c+d), R_{25} = ad/(a+b+c+d), R_{46} = bc/(a+b+c+d),\\
& R_{36} = cd/(a+b+c+d), R_{56} = bd/(a+b+c+d). \end{align*}
\end{enumerate}
\end{theorem}

\begin{proof}
Из примеров~\ref{Box_response} и \ref{Box_drop} знаем ответ и перепад для Box-цепи. Согласно аксиомам $(C), (I), (P), (B),$ система уравнений для H-цепи:
\begin{equation*}
	\begin{cases}
U_1-U_4=\Delta U_{14},\\
U_2-U_3=\Delta U_{23},\\
I_{15}ab/(a+b+c+d)=U_1-U_5,\\
I_{46}bc/(a+b+c+d)=U_4-U_6,\\
I_{56}bd/(a+b+c+d)=U_5-U_6,\\
I_{25}ad/(a+b+c+d)=U_2-U_5,\\
I_{36}cd/(a+b+c+d)=U_3-U_6,\\
I_{15}+I_{46}=0,\\
I_{15}+I_{25}-I_{56}=0,\\
I_{46}+I_{56}+I_{36}=0.
\end{cases}
\end{equation*}
Откуда получаем:
\begin{align*}
&I_{15}=I_{64}=((a+b+c)\Delta{U_{14}}-b\Delta{U_{23}})/(b(a+c)),\\
&I_{25}=I_{63}=((a+c+d)\Delta{U_{23}}-d\Delta{U_{14}})/(d(a+c)),\\
&U_{1}-U_{2}=(\Delta{U_{14}}-\Delta{U_{23}})a/(a+c).
\end{align*}

Перепады и ответы у цепей одинаковы, поэтому они П-эквивалентны.
\end{proof}	

Следующие следствия получаются аналогично доказательству теоремы об электрических преобразованиях \cite[доказательство теоремы об электрических преобразованиях, стр. 35]{SSU}.
\begin{corollary}\label{Box-H transformation1} (рис.~\ref{case2new})
Пусть в плоской двухпортовой цепи есть вершины $x, y, z, t$ и ребра $xy, yz, zt, tx$, причём $I_{xy}+I_{tz}=0$ при любых разностях входных потенциалов.
Тогда если добавить 2 новых неграничных вершины v и w, а также заменить ребра $xy, yz, zt, tx$ на ребра $xv, yv, tw, zw, vw$ сопротивлениями
	
\begin{align}\label{eq1} 
	& R_{xv} = R_{xy}R_{xt}/(R_{xy}+R_{xt}+R_{tz}+R_{yz}),\\
	\nonumber
	&R_{yv} = R_{xy}R_{yz}/(R_{xy}+R_{xt}+R_{tz}+R_{yz}), \\
	\nonumber
	&R_{tw} = R_{xt}R_{tz}/(R_{xy}+R_{xt}+R_{tz}+R_{yz}),\\
	\nonumber
	& R_{zw} = R_{tz}R_{yz}/(R_{xy}+R_{xt}+R_{tz}+R_{yz}),\\ 
	\nonumber
	&R_{vw} = R_{xt}R_{yz}/(R_{xy}+R_{xt}+R_{tz}+R_{yz}), \end{align}
то получится двухпортовая цепь, П-эквивалентная исходной.
\end{corollary}

\begin{corollary}\label{Box-H transformation} (см. первое преобразование на рисунке~\ref{case2new})
Пусть в плоской двухпортовой цепи есть вершины $x, y, z, t$ и ребра $xy, yz, zt, tx$. 
Пусть при удалении ребер $xy, yz, zt, tx$ любые две вершины, где первая --- это одна из $1, 4, x, t$, а вторая --- одна из $2, 3, y, z$, оказываются в разных компонентах связности.
Тогда если добавить 2 новых неграничных вершины $v$ и $w$, а также заменить ребра $xy, yz, zt, tx$ на ребра $xv, yv, tw, zw, vw$ с сопротивлениями, заданными формулой~\eqref{eq1}, то получится двухпортовая цепь, П-эквивалентная исходной.
\end{corollary}

Такое преобразование цепи назовём \emph{Box-H преобразованием}, применённым к рёбрам $xy, yz, zt, tx$.

\begin{theorem}\label{equivalence H or I}
Для любой плоской двухпортовой цепи $G$ найдется П-эквивалентная ей плоская двухпортовая цепь $G'$, удовлетворяющая следующим двум условиям:
\begin{enumerate}[(1)]
	\item[(1)] У цепи $G'$ либо не более 4 ребер, либо это H-цепь.
	\item[(2)] Цепь $G'$ может быть получена из цепи G последовательностью элементарных преобразований (рис.~\ref{elementary}), объединений вершин с одинаковым потенциалом (рис.~\ref{same_potentials2} и предложение~\ref{same potential}) и Box-H преобразований (рис.~\ref{Box-H_transformation}).
\end{enumerate}
\end{theorem}

\begin{remark}\label{rational}
	Если у цепи $G$ все ребра имели рациональные сопротивления, то и у цепи $G'$ сопротивления ребер будут рациональны.
\end{remark}

Если не требовать, чтобы перепад сохранялся, а следить только за ответом, то можно ограничиться лишь 3 ребрами:

\begin{remark}\label{equivalence 3}
Для любой плоской двухпортовой цепи существует плоская двухпортовая цепь с таким же ответом, которая имеет не более 3 ребер. (Можно доказать это, воспользовавшись теоремой~\ref{equivalence H or I} и равенством ответов у цепей, изображенных на рисунке~\ref{N-N}.)
\end{remark}

\begin{proof}[Доказательство теоремы~\ref{equivalence H or I}.] 
	Рассмотрим плоские двухпортовые цепи с наименьшим числом ребер, которые можно получить из цепи $G$ указанными в теореме преобразованиями. Среди таких цепей выберем цепь с наибольшим числом вершин. Тогда в полученной \emph{минимальной} цепи нет циклов длины 3, иначе при помощи преобразования треугольник-звезда можно было бы увеличить число вершин. Кроме того, нет неграничных вершин степени два, иначе можно было бы уменьшить число рёбер последовательным соединением.
	
	По теореме~\ref{Kenyon} и предложению~\ref{equel=>} у этой цепи не более 6 рёбер. Если у этой цепи не более 4 ребер, то теорема доказана. Иначе у этой цепи 5 ребер либо 6. Рассмотрим два случая.

\textit{Случай 1: в цепи 5 ребер.} Эта цепь либо дерево, либо в ней есть циклы длины больше трех.

\textit{Подслучай 1: в цепи нет циклов.} Тогда в цепи 6 вершин. Есть две неграничных вершины 5 и 6 степени не меньше 3. Всего ребер 5, следовательно, есть ребро 56. Тогда это H-цепь или \emph{I-цепь} (цепь на 6 вершинах с 5 ребрами: ${15}, {45}, {26}, {36}, {56},$ рис.~\ref{same_potentials2}). Если это I-цепь, то в силу изолированности портов (условие (P) из \S\ref{prelim2}) $I_{15}=I_{54}$, а значит по правилу Кирхгофа (условие (I)) $I_{56} = 0$. Тогда по закону Ома (условие (C)), $U_5 = U_6$, то есть можно уменьшить число ребер объединением вершин 5 и 6 с одинаковым потенциалом (предложение~\ref{same potential}).
			
\textit{Подслучай 2: в цепи есть цикл.} Цикла длины 5 не может быть, так как тогда была бы неграничная вершина степени два. Значит, есть цикл длины 4 и ещё одно ребро, ведущее в вершину вне цикла (рис.~\ref{case2new}). Неграничной вершиной будет вершина степени 3 (иначе можно уменьшить число ребер). Без ограничения общности, в цикл не входит вершина $4$. Сделаем Box-H преобразование (теорема~\ref{Box-H}), а затем последовательное соединение (предложение~\ref{series transformation}) и получим H-цепь (рис.~\ref{case2new}).

\begin{figure}[h]
	\center
	\includegraphics[width=11cm]{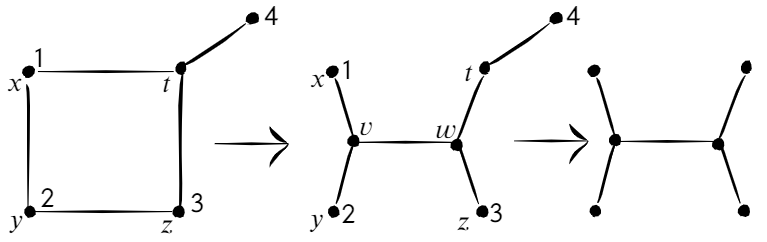}
	\caption{Подслучай 2: цепь с 5 рёбрами и циклом длины 4.}
	\label{case2new}
\end{figure}

\textit{Случай 2: в цепи 6 рёбер.}
По лемме~\ref{>=2 edges} (лемма доказана сразу после доказательства теоремы) степень каждой граничной вершины не меньше 2. У каждой неграничной вершины степень не меньше 3. Значит, эта цепь не дерево, так как нет вершин степени 1. То есть, найдется цикл длины хотя бы 4.

Цикла длины 6 не может быть, так как у неграничных вершин степень хотя бы 3. Цикла длины 5 не может быть, так как шестое ребро будет выходить из неграничной вершины либо в вершину степени один, либо образуя цикл длины 3. 

Если есть цикл длины 4 с ребрами $xy$, $yz$, $zt$, $tx$, то найдётся ребро $xw$, идущее в какую-то вершину вне цикла. Вершин степени 1 и циклов длины 3 нет, значит есть ребро $wz$ (рис.~\ref{xyztw}). Тогда $w$, $y$ и $t$ граничные вершины, так как их степень два. Но тогда цепь не плоская. Это противоречие доказывает теорему (по модулю следующей леммы).
\end{proof}
\begin{figure}[h]
	\center
	\includegraphics[width=4cm]{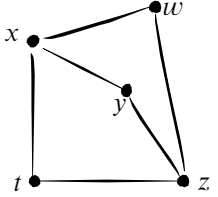}
	\caption{Случай 2: цепь с 6 рёбрами и циклом длины 4.}
	\label{xyztw}
\end{figure}

\begin{lemma}\label{>=2 edges}
	Пусть в двухпортовой цепи 6 рёбер, и из какой-то граничной вершины выходит ровно одно ребро. Тогда можно уменьшить число ребер в цепи указанными в теореме~\ref{equivalence H or I} преобразованиями.
\end{lemma}

\begin{proof}
	Пусть без ограничения общности из вершины 1 выходит ровно одно ребро. Это ребро 12, 13, 14, либо ребро, ведущее в неграничную вершину, которую можно обозначить через 5.
	
	\textit{Случай 1: ребро 12 (рис.~\ref{12}).} Рассмотрим электрическую цепь с 4 входами c ребрами с такими же сопротивлениями как у двухпортовой цепи. Если удалить ребро 12 и вершину 1, то останется цепь с 3 входами и 5 ребрами, которую по теореме~\ref{Kenyon} можно эквивалентно преобразовать в цепь с не более чем 3 ребрами. Если вернуть ребро 12 и вершину 1, то те же преобразования останутся эквивалентными ($I_{12}$ зависит только от потенциалов $U_1$ и $U_2$, а остальные токи от $U_2$, $U_3$, $U_4$). Соответствующая двухпортовая цепь при таких преобразованиях поменяется на П-эквивалентную (предложение~\ref{equel=>}), а количество ребер уменьшится.
	
	\begin{figure}[h!]
		\center
		\includegraphics[width=8cm]{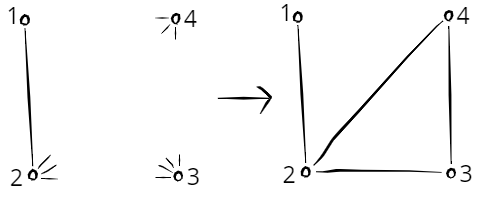}
		\caption{Случай 1: ребро 12.}
		\label{12}
	\end{figure}
	
	\textit{Случай 2: ребро 14.} Аналогичен случаю 1.
	
	\textit{Случай 3: ребро 13 (рис.~\ref{13}).} Вершины плоской цепи расположены на круге, ребро 13 делит его на 2 части -- правую и левую. Рассмотрим электрическую цепь с двумя граничными вершинами 3 и 4, а также ребрами и вершинами такими же, как лежащие правее ребра 13 в двухпортовой цепи. По теореме~\ref{Kenyon} можно преобразовать эту цепь в цепь с одним ребром. Точно так же сделаем с левой частью. Тогда такими преобразованиями получим в цепь, в которой не более 3 ребер.
	
	\begin{figure}[h!]
		\center
		\includegraphics[width=7cm]{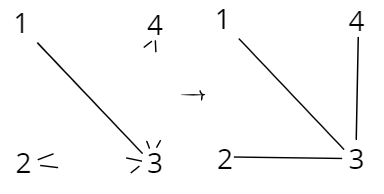}
		\caption{Случай 3: ребро 13.}
		\label{13}
	\end{figure}
	
	\textit{Случай 4: ребро 15 (рис.~\ref{15}).} Рассмотрим часть цепи, которая останется если удалить ребро 15. По теореме~\ref{equivalence H or I} для уже доказанного случая цепей с пятью ребрами эту часть можно П-эквивалентно преобразовать в цепь с не более чем четырьмя ребрами или в H-цепь. Но в последнем случае можно применить последовательное соединение с ребром 15 и уменьшить количество рёбер.
\end{proof}
\begin{figure}[h!]
	\center
	\includegraphics[width=6cm]{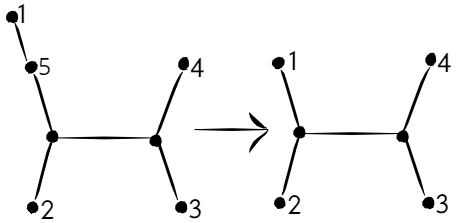}
	\caption{Случай 4: ребро 15.}
	\label{15}
\end{figure}

Теорема~\ref{equivalence H or I} доказана, но возникает вопрос: нельзя ли усилить её и заменить формулировку пункта (1) на такую: \textit{<<У цепи $G'$ не более 4 ребер>>}? Покажем, что нельзя. Более того, покажем, что уменьшить число рёбер в цепи до 4 нельзя никакими преобразованиями, сохраняющими ответ, перепад и <<рациональность сопротивлений>> двухпортовой цепи.

\textit{Вторым перепадом} двухпортовой цепи назовём отображение $\mathds{R}^2 \to \mathds{R}$: $(\Delta{U_{14}}, \Delta{U_{23}}) \mapsto U_{3}-U_{4}$. Для П-эквивалентных цепей вторые перепады равны, так как $U_{3}-U_{4} = \Delta{U_{14}} - \Delta{U_{23}} - (U_{1}-U_{2}).$ \textit{Значением перепада} и \textit{значением второго перепада} назовем числа $U_1-U_2$ и $U_3-U_4$ соответственно.

\begin{predl}\label{contrexample} 
Пусть дана замкнутая H-цепь с ребрами с рациональными сопротивлениями. И пусть оказалось, что отношения любых двух чисел среди значения перепада ($U_1-U_2$), второго перепада ($U_3-U_4$), входных токов ($I_1$ и $I_2$) и разностей входных потенциалов ($\Delta{U_{14}}$ и $\Delta{U_{23}}$) иррациональны.

Тогда любая плоская замкнутая двухпортовая цепь с ребрами с рациональными сопротивлениями, с такими же входными разностями потенциалов, таким же перепадом и таким же ответом, как у Н-цепи, имеет не менее пяти ребер.
\end{predl}

\begin{proof}
Пусть, от противного, существует такая замкнутая двухпортовая цепь $H'$ менее чем с 5 ребрами. Рассмотрим среди цепей, которые можно получить из $H'$ преобразованиями из теоремы~\ref{equivalence H or I}, такую двухпортовую цепь с наименьшим числом ребер, а среди таких --- с наибольшим числом вершин. Тогда в ней три или четыре ребра с рациональными сопротивлениями. Так как для П-эквивалентных цепей перепад, второй перепад и ответ совпадают, то все отношения из условия будут также иррациональны. Рассмотрим два случая в зависимости от наличия циклов в получившейся цепи.

\textit{Случай 1: в цепи нет циклов.} Тогда это дерево на 4 или 5 вершинах. У неграничной вершины не может быть степень один, значит, из какой-то граничной выходит ровно одно ребро. Пусть, без ограничения общности, это ребро выходит из вершины 1. Заметим, что это ребро в вершину 2 или 4 вести не может, так как сопротивление этого ребра рационально, а $(U_1-U_2)/I_1$ и $(U_1-U_4)/I_1$ иррационально. Тогда есть ребро 13 или 15.

\textit{Подслучай 1.1: ребро 13.} Ребро 13 делит плоскую цепь на две части. Без ограничения общности в части с граничной вершиной 4 нет неграничных вершин. В этой части цепи могут быть только ребра 14 или 34, так как цепь плоская. Поскольку в цепи есть ребро 13 и нет циклов, то из вершины 4 выходит единственное ребро 14 или 34 с рациональным сопротивлением. Но $(U_1-U_4)/I_1$ и $(U_3-U_4)/I_1$ иррациональны. Противоречие.

\textit{Подслучай 1.2: ребро 15.} У неграничной вершины 5 степень не меньше 3. Без ограничения общности из неё выходят ребра 15, 25, 35. Тогда степень вершины 4 равна одному и нет ребер 14 и 34, так как сопротивление этих ребер должно быть рационально, а $(U_1-U_4)/I_1$ и $(U_3-U_4)/I_1$ иррациональны. Тогда есть ребро 45 (рис.~\ref{X}). В силу рациональности сопротивлений ребер 15 и 45, а также изолированности портов возникает противоречие с иррациональностью отношения $\Delta{U_{14}}/I_1$.

\begin{figure}[h]
	\center
	\includegraphics[width=4cm]{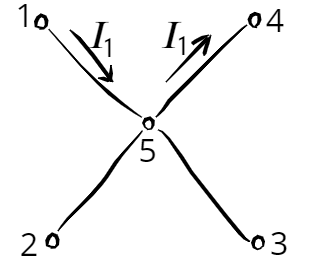}
	\caption{ Подслучай 1.2.}
	\label{X}
\end{figure}

\textit{Случай 2: в цепи есть цикл.} В силу отсутствия циклов длины 3 это Box-цепь. Рациональное сопротивление ребра 12 равно $(U_1-U_2)/(I_1-I_{14})$, рациональное сопротивление ребра 34 равно $(U_3-U_4)/(I_1-I_{14})$. Получаем противоречие с иррациональностью отношения $(U_1-U_2)/(U_3-U_4)$.
\end{proof}

\begin{example}\label{example_min=5}
Рассмотрим H-цепь с сопротивлениями $R_{15}=4, R_{25}=1, R_{46}=4, R_{36}=2, R_{56}=2$.
Любая цепь, которую можно получить из этой цепи при помощи преобразований, указанных в теореме~\ref{equivalence H or I}, имеет не менее пяти ребер, так как при разностях входных потенциалов $\Delta{U_{14}} = 10 - 2 \sqrt{2}$ и $\Delta{U_{14}} = 8 - 5 \sqrt{2}$ будут выполняться условия предложения~\ref{contrexample}.
\end{example}

\section{Разрезания восьмиугольников}\label{tiling_octagons}
Пусть есть 2 прямоугольника $ABCD$ и {$A'B'C'D'$}, причем $A$ и $D$ лежат на отрезке {$A'D'$} и прямоугольники пересекаются только по отрезку {$AD$}. Тогда объединение этих прямоугольников называется \textit{Т-образным восьмиугольником} (рис.~\ref{TZ}).

Пусть есть 2 прямоугольника $ABCD$ и {$A'B'C'D'$}, причем $A$ или $D$ лежит на отрезке {$A'D'$} и прямоугольники пересекаются только по отрезку {$AD'$} или {$A'D$}. Тогда объединение этих прямоугольников называется \textit{Z-образным восьмиугольником} (рис.~\ref{TZ}).

Теорема~\ref{final_theorem} верна не только для П-образного восьмиугольника, но и для Т-образного и Z-образного восьмиугольников. Докажем её в этой, несколько большей, общности с помощью теоремы~\ref{equivalence H or I} и следующей леммы:

\begin{figure}[h]
	\center
	\includegraphics[width=16cm]{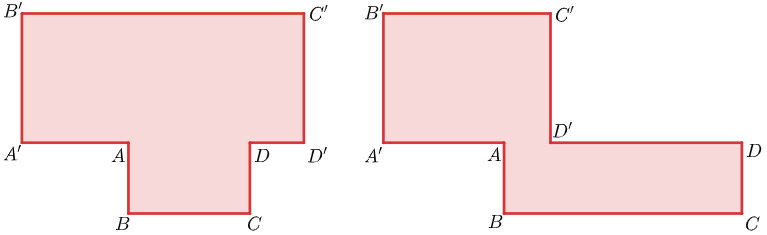}
	\caption{ Т-образный и Z-образный восьмиугольники.}
	\label{TZ}
\end{figure}
\begin{lemma} (Следствие леммы~\ref{correspondence-general}.) \label{correspondence-general2} Пусть $P$ --- ортогональный восьмиугольник общего положения с $4$ вертикальными сторонами, имеющими ориентированные длины ${I_1, I_2, I_3, I_4}$ и x-координаты ${U_1, U_2, U_3, U_4}$. Если $I_1+I_4=0$, то следующие $2$ условия эквивалентны\textup{:}
	\begin{enumerate}[(1)]
		\item \label{1-lemma2} восьмиугольник $P$ можно разрезать на $m$ прямоугольников с отношениями горизонтальной стороны к вертикальной, равными $R_1, \dots, R_m$\textup{;}
		\item \label{2-lemma2} существует плоская замкнутая двухпортовая цепь, состоящая из $m$ существенных рёбер сопротивлениями ${R_1,\dots,R_m}$, с разностями входных потенциалов $\Delta{U_{14}} := U_1-U_4$ и $\Delta{U_{23}} := U_2-U_3$, перепадом $U_1-U_2$ и входными токами ${I_1, I_2}$.
	\end{enumerate}
В частности, все такие двухпортовые цепи для данного восьмиугольника Р являются П-эквивалентными.
\end{lemma}
\begin{proof}
	
$(\ref{1-lemma2})\Rightarrow(\ref{2-lemma2})$. Если выполнено (\ref{1-lemma2}), то по лемме~\ref{correspondence-general} существует плоская замкнутая электрическая цепь с $4$ граничными вершинами, состоящая из $m$ существенных рёбер сопротивлениями ${R_1,\dots,R_m>0}$, с входными потенциалами ${U_1, U_2, U_3, U_4}$ и входными токами ${I_1, I_2, I_3, I_4}$. При этом $I_1+I_4=0$, а значит, эта цепь двухпортовая и выполняется условие (\ref{2-lemma2}).

$(\ref{2-lemma2})\Rightarrow(\ref{1-lemma2})$. Рассмотрим замкнутую плоскую электрическую цепь с входными потенциалами ${U_1, U_2, U_3, U_4}$ и с такими же рёбрами и такими же сопротивлениями, как у двухпортовой цепи из условия (\ref{2-lemma2}). Токи в обеих цепях одинаковые. Тогда по лемме~\ref{correspondence-general} выполняется условие (\ref{1-lemma2}).

\end{proof}

\begin{proof}[Доказательство теоремы~\ref{final_theorem}]
Занумеруем стороны восьмиугольника, начиная с $B'A'$. Тогда для П-, Т- или Z-образного восьмиугольника будет выполняться равенство $I_1+I_4=0$. По лемме~\ref{correspondence-general2}, если такой восьмиугольник можно разрезать на $m$ квадратов, то выполняется (\ref{2-lemma2}).

По теореме~\ref{equivalence H or I} и замечанию~\ref{rational} эта плоская двухпортовая цепь П-эквивалентна цепи с не более чем $5$ рёбрами c рациональными сопротивлениями. В силу П-эквивалентности у новой цепи при разностях входных потенциалов $\Delta{U_{14}}$ и $\Delta{U_{23}}$ будет такой же перепад $U_1-U_2$ и такие же входные токи ${I_1, I_2}$. Тогда по лемме~\ref{correspondence-general2} восьмиугольник можно разрезать на не более 5 прямоугольников с рациональным отношением сторон, так как сопротивления рёбер рациональны.
\end{proof}

В условии теоремы~\ref{final_theorem} число 5 нельзя заменить на меньшее.

\begin{figure}[h]
	\center
	\includegraphics[width=9cm]{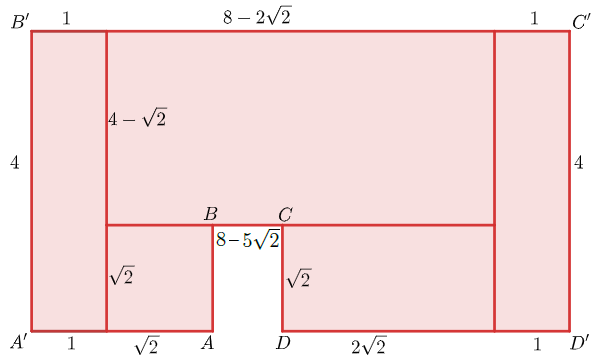}
	\caption{П-образный восьмиугольник, который нельзя разрезать менее чем на 5 прямоугольников с рациональными отношениями сторон.}
	\label{example2}
\end{figure}

\begin{example}\label{example_tiling>4} (Pис.~\ref{example2})
П-образный восьмиугольник со сторонами $A'B'=4, B'C'=10-2\sqrt{2}, AA'=1+\sqrt{2}, AB=\sqrt{2}, DD'=1+2\sqrt{2}$ можно разрезать на 5 прямоугольников с рациональным отношением сторон, но нельзя разрезать менее чем на 5 прямоугольников с рациональным отношением сторон.
\end{example}

\begin{proof}
	Разрезание на 5 прямоугольников показано на рисунке~\ref{example2}. По лемме~\ref{correspondence-general2} существует плоская замкнутая двухпортовая цепь с $5$ ребрами $15, 25, 36, 46,$ и $56$ с рациональными сопротивлениями, такая что $I_1 = 4, I_2 = -\sqrt {2}, I_3 = \sqrt{2}, I_4 = -4, U_1-U_2=-(1+\sqrt{2}), U_3-U_4=-(1+2\sqrt{2}), U_1- U_4=-(1+\sqrt{2}), U_3-U_4=-(1+2\sqrt{2}), \Delta{U_{14}} := U_1-U_4,$ и $\Delta {U_{23}} := U_2-U_3$. Эта цепь удовлетворяет условиям предложения~\ref{contrexample}.
	
	Пусть, от противного, этот П-образный восьмиугольник можно разрезать на менее чем 5 прямоугольников с рациональным отношением сторон. Тогда, аналогично предыдущему абзацу, по лемме~\ref{correspondence-general2} получим двухпортовую цепь с менее чем 5 рёбрами и такими же входными разностями потенциалов, таким же перепадом и таким же ответом. Возникает противоречие с предложением~\ref{contrexample}.
\end{proof}
	
\section{Благодарности}
Автор выражает благодарность научному руководителю, Скопенкову Михаилу Борисовичу, за полезные обсуждения и помощь в написании данной статьи. А также Павлу Пилявскому, который прочитал данную работу и подсказал идеи для дальнейших размышлений.

\bibliographystyle{elsarticle-num}